\documentclass[12pt,oneside,reqno]{amsart}
\usepackage{amsmath}
\usepackage{graphicx}
\usepackage{mathrsfs}
\usepackage{stmaryrd}
\usepackage{amsfonts}
\usepackage{enumerate,amsmath,amssymb,amsthm}
\newcommand{\be}{\begin{eqnarray}}
\newcommand{\ee}{\end{eqnarray}}
\newcommand{\ce}{\begin{eqnarray*}}
\newcommand{\de}{\end{eqnarray*}}
\newtheorem{theorem}{Theorem}[section]
\newtheorem{lemma}[theorem]{Lemma}
\newtheorem{remark}[theorem]{Remark}
\newtheorem{definition}[theorem]{Definition}
\newtheorem{proposition}[theorem]{Proposition}
\newtheorem{Examples}[theorem]{Examples}
\newtheorem{corollary}[theorem]{Corollary}

\def\a{\alpha}

\def\b{\beta}

\def\eps{\epsilon}

\def\[{{\Big[}}
\def\]{{\Big]}}
\def\<{{\langle}}
\def\>{{\rangle}}
\def\({{\Big(}}
\def\){{\Big)}}

\def\bx{{\mathbf{x}}}

\def\dif{{\mathord{{\rm d}}}}

\def\b{\beta}

\def\bt{\begin{theorem}}
\def\et{\end{theorem}}
\def\bl{\begin{lemma}}
\def\el{\end{lemma}}
\def\br{\begin{remark}}
\def\er{\end{remark}}
\def\bx{\begin{Examples}}
\def\ex{\end{Examples}}
\def\bd{\begin{definition}}
\def\ed{\end{definition}}
\def\bp{\begin{proposition}}
\def\ep{\end{proposition}}
\def\bc{\begin{corollary}}
\def\ec{\end{corollary}}
\def\cA{{\mathcal A}}
\def\cB{{\mathcal B}}
\def\cC{{\mathcal C}}
\def\cD{{\mathcal D}}

\def\cF{{\mathcal F}}

\def\cH{{\mathcal H}}

\def\cO{{\mathcal O}}

\def\mB{{\mathbb B}}

\def\mE{{\mathbb E}}

\def\mN{{\mathbb N}}

\def\mR{{\mathbb R}}
\def\mS{{\mathbb S}}

\def\mU{{\mathbb U}}

\def\geq{\geqslant}
\def\leq{\leqslant}

\def\sF{{\mathscr F}}

\def\bE{{\mathbf E}}
\def\bP{{\mathbf P}}

\def\vph{\varphi}
\allowdisplaybreaks

\begin{document}

\title{General Large Deviations and Functional Iterated Logarithm Law for Multivalued Stochastic Differential Equations}

\author{Jiagang Ren$^{1}$, Jing Wu$^{1}$, Hua Zhang$^{2}$}


\date{}

\dedicatory{$^{1}$School of Mathematics and Computational Science,
Sun Yat-Sen University,\\
Guangzhou, Guangdong 510275, P.R.China\\
$^{2}$School of Statistics, Jiangxi University of Finance and Economics,\\
Nanchang, Jiangxi 330013, P.R.China\\
Emails: J. Ren: renjg@mail.sysu.edu.cn\\
J. Wu: wjjosie@hotmail.com\\
H. Zhang: zh860801@163.com}


\thanks{This work is supported by NSFs of China (No. 11171358, 11101441 and 11301553), Doctor Fund of Ministry of Education (No. 20100171110038, 20100171120041 and 20120171120008) and and China Postdoctoral Science Foundation (No. 2013T60817).}

\begin{abstract}
In this paper, we prove a large deviation principle of Freidlin-Wentzell's type for the
multivalued stochastic differential equations, which is a little more general than the results obatined by Ren, Xu and Zhang \cite{RXZ}. As an application, we derive a functional iterated logarithm law for the solutions of multivalued stochastic differential equations.
\end{abstract}

\maketitle

\section{Introduction}

In \cite{RXZ} a large deviation principle of Freidlin-Wentzell type is established for the following multivalued stochastic differential equation (MSDE):
\be\label{eq0} \left\{
\begin{array}{ll}
\dif X(t)\in b(X(t))\dif t+\sigma(X(t))\dif W(t)-A(X(t))\dif t,\\
X(0)=x\in \overline{D(A)},
\end{array}
\right.
\ee
where $A$ is a multivalued maximal monotone operator with domain $D(A):=\{x\in \mR^m:A(x)\neq \emptyset\}$ and graph
$\mathrm{Gr}(A):=\{(x,y)\in \mR^{2m}: x\in \mR^m, y\in A(x)\}$,
$W(t)=\{W^k(t), t\geq 0, k\in\mN\}$ is a sequence of independent standard Brownian motions
on a filtered probability space $(\Omega,\cF,\bP; (\cF_t)_{t\geq 0})$,
$b: \mR^m\to\mR^m$ and $\sigma:\mR^m\to\mR^m\times l^2$ are two continuous functions, $l^2$ stands for the Hilbert space of square summable sequences of real numbers.

We do not recall here the notions of multivalued maximal monotone operators and multivalued stochastic differential equations etc, but instead
refer the readers to \cite{C1, C2, RXZ} for them.

In the present paper, inspired by the work \cite{BC}, we consider the small perturbation of the following type:
\be\label{pmsde}
\left\{
\begin{array}{ll}
\dif X^\eps(t)\in b_\eps(X^\eps(t))\dif t+\sqrt{\eps}
\sigma_\eps(X^\eps(t))\dif W(t)-A_\eps(X^\eps(t))\dif t,\\
X^\eps(0)=x\in \overline{D(A)},\ \ \eps\in(0,1].
\end{array}
\right.
\ee

Compared with \cite{BC}, its novelty  is the appearance of $A$ and $A_{\epsilon}$. Of course, for a large deviation principle to hold for such a problem it is necessary that all the coefficients should converge in some way and thus one has to know how to measure the difference between them. This poses no problem concerning $\sigma_\eps$ and $b_\eps$ as is done in \cite{BC}, but brings some troubles for the singular unbounded $A_\eps$ and $A$.  The point is that they are all multivalued (i.e., set-valued) operators and there seems no straightforward way to measure the difference between them. To circumvent this difficulty we turn to measuring their Yosida approximations by assuming that
\ce
\lim_{\epsilon\rightarrow0}(1+\alpha A_{\epsilon})^{-1}y=(1+\alpha A)^{-1}y,\quad\forall \alpha>0
\de
uniformly for $y$ in compact sets. Our main result will be proved under this condition.

Such an assumption is quite reasonable since in the case of subdifferentials of convex functions it is actually implied by the point-wise convergence of the associated convex functions, see Example 4.4 below. But the price we have to pay is that we have to use the Yosida approximating equations in order to get control of them, which, nonetheless, do not even appear in \cite{RXZ}.

We now give an outline of the paper. In Section 2, we recall the well-known Laplace principle. In Section 3, we present our main result, and we put a detailed proof in Section 4. As an application, in Section 5 we derive the functional iterated logarithm law, which has been obtained through a different method in \cite{Rab} when $A=\partial I_O$, where $O=\{x\in\mR^m,x=(x^1,\cdots,x^m),x^1\geq0\}$.

\section{Preliminaries}
First we recall basic some notions and results. By a solution of (\ref{eq0}) we mean a pair $(X,K)$ of $(\sF_t)$-adapted continuous processes satisfying
\begin{enumerate}[(i)]
\item  $X(0)=x$ and for all $t\geq0$, $X(t)\in \overline{D(A)}$ a.s.;
\item $K$ is of locally finite variation and $K(0)=0$ a.s.;
\item $\dif X(t)=b(X(t))\dif t+\sigma(X(t))\dif W(t)-\dif K(t)$,
$0\leq t<\infty,
\quad a.s.$;
\item For any continuous and $(\sF_t)-$adapted functions $(\a,\b)$
with $(\a(t),\b(t))\in \mathrm{Gr}(A),~ \forall t\in[0,+\infty)$,
 $\<X(t)-\a(t), \dif K(t)-\b(t)\dif t\>\geq0,\quad a.s.$,
 \end{enumerate}
 where $\<\cdot,\cdot\>$ denotes the scalar product in $\mR^d$.

 We will need the following result which is taken from \cite[Proposition 4.1,Proposition 4.4]{C2}
(see also \cite[Propositions 3.3 and 3.4]{Z1}).
\bp\label{bp}
Suppose that $\mathrm{Int}(D(A))\neq\emptyset$. Then for any $a\in\mathrm{Int}(D(A))$, there exist $\gamma>0$, $\mu \geq 0$ such that
for any pair of processes $(X,K)$ solving Eq. (\ref{eq0}) and all $0\leq s<t\leq T$,
\ce
\int_s^t\<X(r)-a,\dif K(r)\>_{\mR^m}\geq \gamma|K|_t^s
-\mu\int_s^t|X(r)-a|\dif r-\gamma\mu(t-s),
\de
where $|K|_t^s$ denotes the total variation of $K$ on $[s,t]$.

Moreover, for any pair of $(X,K)$ and
$(\widetilde X,\widetilde K)$ solving Eq. (\ref{eq0}) with different initial points,
\ce
\<X(t)-\widetilde X(t),\dif K(t)-\dif\widetilde K(t)\>_{\mR^m}\geq0,\quad\text{a.s.}
\de
\ep

We next recall an abstract criterion for Laplace principle, which is equivalent to the large deviation principle (cf. \cite{BD1}). Let $\mU$ be a real separable Hilbert space. It is well known that there exists a Hilbert space so that $l^2\subset\mU$ is Hilbert-Schmidt with embedding operator $J$ and $\{W^k(t),k\in\mN\}$ is a  Brownian motion with values in $\mU$, whose covariance operator is given by $Q=J\circ J^*$. For example, one can take $\mU$ as the completion of $l^2$ with respect to the norm generated by scalar product
\ce
\<h,h'\>_\mU:=\left(\sum_{k=1}^\infty \frac{h_k h'_k}{k^2}\right)^{\frac{1}{2}}, \ \ h,h'\in l^2.
\de

For a Polish space $\mB$, we denote by $\cB(\mB)$ the Borel $\sigma$-field, and by $\mathcal {C}_T(\mB)$ the continuous function space from $[0,T]$
to $\mB$, and endow it with the uniform distance so that $\mathcal {C}_T(\mB)$ is still a Polish space. Define
\ce
\ell^2_T:=\left\{h=\int^\cdot_0\dot h(s)\dif s: ~~\dot h\in L^2([0,T];l^2)\right\}
\de
with the norm
\ce
\|h\|_{\ell^2_T}:=\left(\int^T_0\|\dot h(s)\|_{l^2}^2\dif s\right)^{1/2},
\de
where the dot denotes the generalized derivative. Let $\mu$ be the law of the Brownian motion $W$ in $\mathcal {C}_T(\mU)$. Then $(\ell^2_T,\mathcal{C}_T(\mU),\mu)$ forms an abstract Wiener space.

For $T,N>0$,  set
\ce
\cD_N:=\{h\in \ell^2_T: \|h\|_{\ell^2_T}\leq N\}
\de
and
\ce
\cA^T_N:=\left\{
\begin{aligned}
&\mbox{ $h: [0,T]\to l^2$ is a continuous and
$(\cF_t)$-adapted }\\
&\mbox{ process, and for almost all $\omega$},\ \ h(\cdot,\omega)\in\cD_N
\end{aligned}
\right\}.
\de
We equip $\cD_N$ with the  weak convergence topology in $\ell^2_T$. Then
\ce
\mbox{$\cD_N$ is metrizable as a compact Polish space}.
\de

\bd
Let $\mS$ be a Polish space. A function $I:\mS\rightarrow[0,\infty]$ is called a rate function if for every $a<\infty$, the set
$\{f\in\mS: I(f)\leq a\}$ is compact in $\mS$.
\ed

Let $\{Z^\eps: \mathcal {C}_T(\mU)\to\mS,\eps\in(0,1)\}$ be a family of measurable mappings. Assume that there is a measurable map $Z_0:\ell^2_T\mapsto \mS$ such that
\begin{enumerate}[{\bf (LD)$_\mathbf{1}$}]
\item
For any $N>0$, if a family $\{h_\eps, \eps\in(0,1)\}\subset\cA^T_N$ (as random variables in $\cD_N$) converges in distribution to $h\in \cA^T_N$, then for some subsequence $\eps_k$, $Z^{\eps_k}\Big(\cdot+\frac{h_{\eps_k}(\cdot)}{\sqrt{\eps_k}}\Big)$ converges in distribution to $Z_0(h)$ in $\mS$.
\end{enumerate}
\begin{enumerate}[{\bf (LD)$_\mathbf{2}$}]
\item
For any $N>0$, if $\{h_n,n\in\mN\}\subset \cD_N$ weakly converges to $h\in\ell^2_T$, then for some subsequence $h_{n_k}$, $Z_0(h_{n_k})$ converges to $Z_0(h)$ in $\mS$.
\end{enumerate}

For each $f\in\mS$, we define
\be\label{La}
I(f):=\frac{1}{2}\inf_{\{h\in\ell^2_T;f=Z_0(h)\}}\|h\|^2_{\ell^2_T},\label{ra}
\ee
where $\inf\emptyset=\infty$ by convention. Then under {\bf(LD)$_\mathbf{2}$}, $I(f)$ is a  rate function.

We recall the following result due to \cite{BD1}.
\bt\label{Th2}
Under {\bf (LD)$_\mathbf{1}$} and {\bf (LD)$_\mathbf{2}$}, $\{Z^\eps,\eps\in(0,1)\}$ satisfies the Laplace principle with the rate function $I(f)$ given by (\ref{ra}). More precisely, for each real bounded continuous function $g$ on $\mS$:
\ce
\lim_{\eps\rightarrow 0}\eps\log\bE^{\mu}[\exp\{-\frac{g(Z^\eps)}{\eps}\}]
=-\inf_{f\in\mS}\{g(f)+I(f)\}.
\de
In particular, the family of $\{Z^\eps,\eps\in(0,1)\}$ satisfies  the large deviation principle in $(\mS,\cB(\mS))$ with the rate function $I(f)$. More precisely, let $\nu_\eps$ be the law of $Z^\eps$ in $(\mS,\cB(\mS))$, then for any $B \in\cB(\mS)$
\ce
-\inf_{f\in B^o}I(f)\leq\liminf_{\eps\rightarrow 0}\eps\log\nu_\eps(B)
\leq\limsup_{\eps\rightarrow 0}\eps\log\nu_\eps(B)\leq -\inf_{f\in \bar B}I(f),
\de
where the closure and the interior are taken in $\mS$, and $I(f)$ is defined in (\ref{ra}).
\et

For $0<\eta<1/e$, we define a concave function as
\ce
\rho_{\eta}(x):=
\begin{cases}
x\log x^{-1},&x\leq \eta;\\
\eta\log\eta^{-1}+(\log\eta^{-1}-1)(x-\eta),&x>\eta.
\end{cases}
\de

The following generalization of the Gronwall-Belmman type inequality (cf. \cite{RZ1}) is useful in our paper.
\bl[Bihari's inequality]
If $g(s)$, $q(s)$ are two strictly positive functions on $\mR^+$ such that
\ce
g(t)\leq g(0)+\int_0^tq(s)\rho_{\eta}(g(s))\dif s,\quad t\geq 0.
\de
Then
\ce
g(t)\leq (g(0))^{\exp\{-\int_0^tq(s)ds\}}.
\de
\el

The following lemma can been also found in \cite{RZ1}.
\bl\label{elementary}
\begin{enumerate}
\item $\rho_{\eta}$ is decreasing in $\eta$, i.e., $\rho_{\eta_1}\leq\rho_{\eta_2}$, $\eta_2<\eta_1<1/e$.
\item For any $p>1$ and $\eta$ sufficiently small, we have
\ce
x^p\rho_{\eta}(x)\leq\frac{1}{1+p}\rho_{\eta}(x^{1+p}).
\de
\end{enumerate}
\el
Throughout the paper, $C$ will denote different constants (dependent on the indexes or not) whose values are not important.

\section{Main Results}

We assume that
\begin{enumerate}[{\bf (H1)}]
\item
$A,A_{\epsilon}$ are maximal monotone operators with a common domain $D=\overline{D(A)}=\overline{D(A_{\epsilon})}$ and  $0\in\mathrm{Int}(D(A))$. Moreover, $A_{\epsilon}$ is locally bounded at $0$ uniformly in $\epsilon$, i.e., there exists $\gamma>0$ such that
\ce
\sup\{|y|;y\in A_{\epsilon}(x),\epsilon\in[0,1],x\in B_0(\gamma):=\{x\in\mR^m;|x|\leq\gamma\}\}<\infty,~\forall\epsilon;
\de
\item
$\sigma$ and $\sigma_{\epsilon}$ are continuous functions
and satisfy that for some $C$ and all $x,y\in\mR^m$,
\ce
\|\sigma(x)-\sigma(y)\|_{L_2(l^2;\mR^m)}^2\vee
\|\sigma_{\epsilon}(x)-\sigma_{\epsilon}(y)\|_{L_2(l^2;\mR^m)}^2&\leq& C|x-y|^2(1\vee\log|x-y|^{-1}),\\
\|\sigma_{\epsilon}(x)\|_{L_2(l^2;\mR^m)}\vee\|\sigma(x)\|_{L_2(l^2;\mR^m)}&\leq& C(1+|x|),
\de
where $C$ is independent of $\epsilon$ and $L_2(l^2;\mR^m)$ denotes the Hilbert space of Hilbert-Schmidt operator from $l^2$ to $\mR^m$ and $|\cdot|$ denotes the norm in $\mR^m$, and
\ce
\lim_{\epsilon\rightarrow 0}\|\sigma_{\epsilon}(y)-\sigma(y)\|_{L_2(l^2;\mR^m)}^2=0
\de
uniformly for $y$ in compact sets;
\item
$b$ and $b_{\epsilon}$ are continuous functions and satisfy that for some $C>0$ and all $x,y\in\mR^m$,
\ce
\<x-y,b(x)-b(y)\>\vee
\<x-y,b_{\epsilon}(x)-b_{\epsilon}(y)\>&\leq&C|x-y|^2(1\vee\log|x-y|^{-1}),\\
|b_{\epsilon}(x)|\vee|b(x)|&\leq& C(1+|x|),
\de
where $C$ is independent of $\epsilon$, and
\ce
\lim_{\epsilon\rightarrow 0}|b_{\epsilon}(y)-b(y)|=0
\de
uniformly for $y$ in compact sets;
\item
for all $\alpha>0$,
\ce
\lim_{\epsilon\rightarrow0}(1+\alpha A_{\epsilon})^{-1}y=(1+\alpha A)^{-1}y
\de
uniformly for $y$ in compact sets.
\end{enumerate}

\br
It is obvious that {\bf (H4)} is equivalent to
\ce
\lim_{\epsilon\downarrow 0}A^\alpha_\epsilon(y)=A^\alpha(y)
\de
uniformly for $y$ in compacts sets, where $A^\alpha$ (resp. $A_\eps^\alpha$) is the Yosida approximation of $A$ (resp. $A_\eps$).
\er
\br
By shifting, {\bf (H1)} can be replaced by the existence $V\subset\mathrm{Int}D(A_{\epsilon})$ such that
\ce
\sup\{|y|;y\in A_{\epsilon}(x),\epsilon\in[0,1],x\in V\}<\infty.
\de
\er

It is well known that under {\bf (H1)}-{\bf (H3)}, for every $\epsilon>0$, there exists a unique solution
$(X^\eps(\cdot,x),K^\eps(\cdot,x))$ to Eq. (\ref{pmsde}) (cf. \cite{Z1}). Our main result is stated as follows.
\bt\label{Main}
Assume that {\bf (H1)}-{\bf (H4)} hold.
Then the family of $\{X^\eps(t,x),\ \eps\in(0,1)\}$ satisfies the large deviation principle in $\mS:=\mathcal{C}([0,T]\times\overline{D(A)};\overline{D(A)})$ with the rate function given by
\ce
I(f)=\frac{1}{2}\inf_{\{h\in\ell^2_T;f=X^h\}}\|h\|^2_{\ell^2_T},
\de
where $X^h(t,x)$ solves the following MDE:
\ce
\begin{cases}
\dif X^h(t)\in b(X^h(t))\dif t+
\sigma(X^h(t))\dot h(t)\dif t-A(X^h(t))\dif t,\\
X^h(0)=x.
\end{cases}
\de
More precisely, for any $B\in\cB(\mS)$,
\ce
-\inf_{f\in B^o}I(f)\leq\liminf_{\epsilon\rightarrow0}\epsilon\log \bP\{X^{\epsilon}\in B\}\leq\limsup_{\epsilon\rightarrow0}\epsilon \log \bP\{X^{\epsilon}\in B\}\leq-\inf_{f\in\bar{B}}I(f).
\de
\et

\section{Proof of Theorem \ref{Main}}

This section is devoted to the proof of Theorem \ref{Main}. To prove this result, by Theorem \ref{Th2}, the main task is to verify {\bf(LD)$_\mathbf{1}$} and {\bf (LD)$_\mathbf{2}$} with
\ce
\mS:=\cC([0,T]\times \overline{D(A)};\overline{D(A)}),\ \ Z^\eps=X^\eps,\ \ Z_0(h)=X^h.
\de

In the rest of this section, we suppose that $h_\eps\in\cA^T_N$ converge almost surely to $h\in\cA^T_N$ as random variables in $\ell^2_T$.

\subsection{Yosida approximations}

Consider the following controlled MSDE:
\be\label{control}
\begin{cases}
\dif X^{\epsilon,h_{\epsilon}}(t)\in b_{\epsilon}(X^{\epsilon,h_{\epsilon}}(t))\dif t+\sigma_{\epsilon}(X^{\epsilon,h_{\epsilon}}(t))\dot{h}_{\epsilon}(t)\dif t+\sqrt{\epsilon}\sigma_{\epsilon}(X^{\epsilon,h_{\epsilon}}(t))\dif W(t)-A_{\epsilon}(X^{\epsilon,h_{\epsilon}}(t))\dif t,\\
X^{\epsilon,h_{\epsilon}}(0)=x\in \overline{D(A)}.
\end{cases}
\ee
By Girsanov's transformation result, it admits a unique solution $(X^{\epsilon,h_{\epsilon}}(\cdot,x),K^{\epsilon,h_{\epsilon}}(\cdot,x))$.

Now we consider the following Yosida approximation:
\be\label{Yosida1}
\begin{cases}
\dif X^{\epsilon,h_{\epsilon},\alpha}(t)= b_{\epsilon}(X^{\epsilon,h_{\epsilon},\alpha}(t))\dif t+\sigma_{\epsilon}(X^{\epsilon,h_{\epsilon},\alpha}(t))\dot{h}_{\epsilon}(t)\dif t\\
\quad\quad\quad\quad\quad\quad{}+\sqrt{\epsilon}\sigma_{\epsilon}(X^{\epsilon,h_{\epsilon},\alpha}(t))\dif W(t)-A_{\epsilon}^{\alpha}(X^{\epsilon,h_{\epsilon},\alpha}(t))\dif t,\\
X^{\epsilon,h_{\epsilon},\alpha}(0)=x\in \overline{D(A)}.
\end{cases}
\ee

We have the following lemmas which will play an important role in the proof of Theorem \ref{Main}.
\bl\label{b1}
Assume that {\bf (H1)}-{\bf (H4)} hold. For any $p\geq 1$ and $x\in \overline{D(A)}$, there exists a constant $C>0$ such that for any $\epsilon\in(0,1)$ and $\alpha>0$,
\ce
\bE\sup_{t\in[0,T]}|X^{\epsilon,h_{\epsilon},\alpha}(t,x)|^{2p}\leq C,
\de
where $C$ may be dependent of $p$, $T$, $N$ and $x$, but is independent of $\epsilon$ and $\alpha$.
\el
\begin{proof}
For notational simplicity we omit the index $x$. Since $A_{\epsilon}^{\alpha}$ is monotone, it follows that
\ce
&&-\<X^{\epsilon,h_{\epsilon},\alpha}(s), A_{\epsilon}^{\alpha}(X^{\epsilon,h_{\epsilon},\alpha}(s))\>_{\mR^m}\\
&=&-\<X^{\epsilon,h_{\epsilon},\alpha}(s), A_{\epsilon}^{\alpha}(X^{\epsilon,h_{\epsilon},\alpha}(s))-A_{\epsilon}^{\alpha}(0)
+A_{\epsilon}^{\alpha}(0)\>_{\mR^m}\\
&\leq&-\<X^{\epsilon,h_{\epsilon},\alpha}(s), A_{\epsilon}^{\alpha}(0)\>_{\mR^m}.
\de
Noting that $\sup_\epsilon|A_{\epsilon}^{\alpha}(0)|\leq\sup_\epsilon|A_{\epsilon}^0(0)|$,
\ce
-\<X^{\epsilon,h_{\epsilon},\alpha}(s), A_{\epsilon}^{\alpha}(X^{\epsilon,h_{\epsilon},\alpha}(s))\>_{\mR^m}\leq\frac{1}{2}|X^{\epsilon,h_{\epsilon},\alpha}(s)|^2+C.
\de
Then by It\^o's formula, {\bf (H1)}-{\bf (H3)}, we have
\ce
|X^{\epsilon,h_{\epsilon},\alpha}(t)|^{2p}
&=&|x|^{2p}+2p\int_0^t|X^{\epsilon,h_{\epsilon},\alpha}(s)|^{2p-2}
\<X^{\epsilon,h_{\epsilon},\alpha}(s),b(X^{\epsilon,h_\epsilon,\alpha}(s))\>_{\mR^m}\dif s\\
&&{}+2p\int_0^t|X^{\epsilon,h_{\epsilon},\alpha}(s)|^{2p-2}
\<X^{\epsilon,h_{\epsilon},\alpha}(s),
\sigma_{\epsilon}(X^{\epsilon,h_{\epsilon},\alpha}(s))\dot{h}_{\epsilon}(s)\>_{\mR^m}\dif s\\
&&{}+2p\sqrt{\epsilon}\int_0^t|X^{\epsilon,h_{\epsilon},\alpha}(s)|^{2p-2}
\<X^{\epsilon,h_{\epsilon},\alpha}(s),\sigma_{\epsilon}(X^{\epsilon,h_{\epsilon},\alpha}(s))\dif W(s)\>_{\mR^m}\\
&&{}+p\epsilon\int_0^t|X^{\epsilon,h_{\epsilon},\alpha}(s)|^{2p-2}
(\|\sigma_{\epsilon}(X^{\epsilon,h_{\epsilon},\alpha}(s))\|_{L_2(l^2;\mR^m)}^2\\
&&{}+2(p-1)\frac{\<X^{\epsilon,h_{\epsilon},\alpha}(s),
\sigma_{\epsilon}(X^{\epsilon,h_{\epsilon},\alpha}(s))\sigma_{\epsilon}^*
(X^{\epsilon,h_{\epsilon},\alpha}(s))(X^{\epsilon,h_{\epsilon},\alpha}(t))
\>_{\mR^m}}{|X^{\epsilon,h_{\epsilon},\alpha}(s)|^2})\dif s\\
&&{}-2p\int_0^t|X^{\epsilon,h_{\epsilon},\alpha}(s)|^{2p-2}
\<X^{\epsilon,h_{\epsilon},\alpha}(s),A_{\epsilon}^{\alpha}(X^{\epsilon,h_{\epsilon},
\alpha}(s))\>_{\mR^m}\dif s\\
&\leq&C+C\int_0^t|X^{\epsilon,h_{\epsilon},\alpha}(s)|^{2p}\dif s
\\
&&{}+2p\int_0^t|X^{\epsilon,h_{\epsilon},\alpha}(s)|^{2p-1}|\sigma_{\epsilon}(X^{\epsilon,h_{\epsilon},\alpha}(s))
\dot{h}_{\epsilon}(s)|\dif s\\
&&{}+C\int_0^t|X^{\epsilon,h_{\epsilon},\alpha}(s)|^{2p}\|\dot{h}_\eps(s)\|_{l^2}\dif s\\
&&{}+2p\sqrt{\epsilon}\int_0^t|X^{\epsilon,h_{\epsilon},\alpha}(s)|^{2p-2}
\<X^{\epsilon,h_{\epsilon},\alpha}(s),\sigma_{\epsilon}(X^{\epsilon,h_{\epsilon},\alpha}(s))\dif W(s)\>_{\mR^m}.
\de
Set
$$
f(t):=\bE\sup_{s\in[0,t]}|X^{\epsilon,h_{\epsilon},\alpha}(s)|^{2p}.
$$
By BDG's inequality, {\bf (H2)} and Young's inequality, we have
\ce
&&\bE\sup_{t'\in[0,t]}\int_0^{t'}|X^{\epsilon,h_{\epsilon},\alpha}(s)|^{2p-2}
\<X^{\epsilon,h_{\epsilon},\alpha}(s),\sigma_{\epsilon}(X^{\epsilon,h_{\epsilon},\alpha}(s))\dif W(s)\>_{\mR^m}\nonumber\\
&\leq& C\bE[(\int_0^{t}|X^{\epsilon,h_{\epsilon},\alpha}(s)|^{4p-2}
\|\sigma_{\epsilon}(X^{\epsilon,h_{\epsilon},\alpha}(s)\|_{l^2}^2\dif s)^{1/2}]\nonumber\\
&\leq& C\bE[(\sup_{s\in[0,t]}|X^{\epsilon,h_{\epsilon},\alpha}(s)|^{2p}
\int_0^{t}|X^{\epsilon,h_{\epsilon},\alpha}(s)|^{2p}
\dif s)^{1/2}]+C\int_0^{t}\bE|X^{\epsilon,h_{\epsilon},\alpha}(s)|^{2p}\dif s+C\nonumber\\
&\leq& \frac{1}{8}f(t)+C\int_0^{t}\bE|X^{\epsilon,h_{\epsilon},\alpha}(s)|^{2p}\dif s+C.
\de
Similarly, we have
\ce
\bE[\int_0^t|X^{\epsilon,h_{\epsilon},\alpha}(s)|^{2p}\|\dot{h}_\eps(s)\|_{l^2}\dif s]
&\leq& N\bE[(\int_0^t|X^{\epsilon,h_{\epsilon},\alpha}(s)|^{4p}\dif s)^{1/2}]\nonumber\\
&\leq&\frac{1}{4}f(t)+C_{N}\int_0^{t}\bE|X^{\epsilon,h_{\epsilon},\alpha}(s)|^{2p}\dif s
\de
and
\ce
\bE[\int_0^t|X^{\epsilon,h_{\epsilon},\alpha}(s)|^{2p-1}|\sigma_{\epsilon}(X^{\epsilon,h_{\epsilon},\alpha}(s))
\dot{h}_{\epsilon}(s)|\dif s]\leq\frac{1}{8}f(t)+C\int_0^{t}\bE|X^{\epsilon,h_{\epsilon},\alpha}(s)|^{2p}\dif s+C.
\de
Combining these we get
\ce
f(t)\leq \frac{1}{2}f(t)+C+C\int^t_0\bE|X^{\epsilon,h_{\epsilon},\alpha}(s)|^{2p}\dif s,
\de
and hence
\ce
f(t)\leq 2C+2C\int^t_0f(s)\dif s.
\de
Noting that the constant $C$ is independent of $\epsilon$ and $\alpha$, due to {\bf (H1)}-{\bf (H4)}, we obtain the desired estimate by Gronwall's inequality.
\end{proof}

The following result can be proved in a similar way as above by using Proposition \ref{bp} (see also \cite[Lemma 3.5]{RXZ}).

\bl\label{b2}
Assume that {\bf (H1)}-{\bf (H4)} hold. For any $p\geq 1$ and $x\in \overline{D(A)}$, there exists a constant $C>0$ such that for any $\epsilon\in(0,1)$,
\ce
\bE\sup_{t\in[0,T]}|X^{\epsilon,h_{\epsilon}}(t,x)|^{2p}+\bE|K^{\epsilon,h_{\epsilon}}(\cdot,x)|_T^0\leq C,
\de
where $C$ may be dependent of $p$, $T$, $N$ and $x$, but is independent of $\epsilon$.
\el

\bl\label{p1}
Assume that {\bf (H1)}-{\bf (H4)} hold. Then, for any $x\in \overline{D(A)}$, we have
\ce
\lim_{\alpha\rightarrow0}\bE\sup_{t\in[0,T]}|X^{\epsilon,h_{\epsilon},\alpha}(t,x)-X^{\epsilon,h_{\epsilon}}(t,x)|^2=0
\de
uniformly in $\epsilon>0$.
\el

\begin{proof}
The proof will be carried out by adapting the ideas from \cite{C2,Ras}.

Let
\ce
M_{\epsilon}(t)=\int_0^tb_{\epsilon}(X^{\epsilon,h_{\epsilon}}(s))\dif s+\int_0^t\sigma_{\epsilon}(X^{\epsilon,h_{\epsilon}}(s))\dot{h}_{\epsilon}(s)\dif s+\sqrt{\epsilon}\int_0^t\sigma_{\epsilon}(X^{\epsilon,h_{\epsilon}}(s))\dif W(s).
\de
By Lemma \ref{b2}, BDG's inequality and {\bf (H2)}-{\bf (H3)} we obtain
\be\label{estimate1}
\bE\sup_{t\in[0,T]}|M_{\epsilon}(t)|^2<\infty,
\ee
and
\be\label{estimate2}
\bE\sup_{0\leq|t'-t|\leq\delta}|M_{\epsilon}(t')-M_{\epsilon}(t)|^p\leq C_p\delta^{p/2-1},\quad\forall p\geq 2,
\ee
where $C_p$ is a constant independent of $\delta$ and $\epsilon$.

Now we consider the following Yosida approximation:
\be\label{Yosida2}
\dif \widetilde{X}^{\epsilon,h_{\epsilon},\alpha}(t)=\dif M_{\epsilon}(t)-A_{\epsilon}^{\alpha}(\widetilde{X}^{\epsilon,h_{\epsilon},\alpha}(t))\dif t,\quad
\widetilde{X}^{\epsilon,h_{\epsilon},\alpha}(0)=x.
\ee
Since $A_{\epsilon}^{\alpha}$ is Lipschitz, Eq. (\ref{Yosida2}) admits a unique solution.

First by adapting the arguments in \cite{Z}, we get (see Appendix for a proof)
\be\label{key1}
\lim_{\alpha\rightarrow0}\bE\sup_{t\in[0,T]}|\widetilde{X}^{\epsilon,h_{\epsilon,\alpha}}(t)-X^{\epsilon,h_{\epsilon}}(t)|^2=0
\ee
uniformly in $\epsilon>0$.

Secondly, by It\^{o}'s formula and the monotonicity of $A_{\epsilon}$, we have
\ce
&&|\widetilde{X}^{\epsilon,h_{\epsilon},\alpha}(t)-X^{\epsilon,h_{\epsilon},\alpha}(t)|^2\\
&\leq&2\int_0^t\<\widetilde{X}^{\epsilon,h_{\epsilon},\alpha}(s)-X^{\epsilon,h_{\epsilon}}(s),b_{\epsilon}(X^{\epsilon,h_{\epsilon}}(s))-b_{\epsilon}(X^{\epsilon,h_{\epsilon},\alpha}(s))\>_{\mR^m}\dif s\\
&&{}+2\int_0^t\<X^{\epsilon,h_{\epsilon}}(s)-X^{\epsilon,h_{\epsilon},\alpha}(s),b_{\epsilon}(X^{\epsilon,h_{\epsilon}}(s))-b_{\epsilon}(X^{\epsilon,h_{\epsilon},\alpha}(s))\>_{\mR^m}\dif s\\
&&{}+2\int_0^t\<\widetilde{X}^{\epsilon,h_{\epsilon},\alpha}(s)-X^{\epsilon,h_{\epsilon}}(s),\sigma_{\epsilon}(X^{\epsilon,h_{\epsilon}}(s))\dot{h}_{\epsilon}(s)-\sigma_{\epsilon}(X^{\epsilon,h_{\epsilon},\alpha}(s))\dot{h}_{\epsilon}(s)\>_{\mR^m}\dif s\\
&&{}+2\int_0^t\<X^{\epsilon,h_{\epsilon}}(s)-X^{\epsilon,h_{\epsilon},\alpha}(s),\sigma_{\epsilon}(X^{\epsilon,h_{\epsilon}}(s))\dot{h}_{\epsilon}(s)-\sigma_{\epsilon}(X^{\epsilon,h_{\epsilon},\alpha}(s))\dot{h}_{\epsilon}(s)\>_{\mR^m}\dif s\\
&&{}+2\sqrt{\epsilon}\int_0^t\<\widetilde{X}^{\epsilon,h_{\epsilon},\alpha}(s)-X^{\epsilon,h_{\epsilon}}(s),(\sigma_{\epsilon}(X^{\epsilon,h_{\epsilon}}(s))-\sigma_{\epsilon}(X^{\epsilon,h_{\epsilon},\alpha}(s)))\dif W(s)\>_{\mR^m}\\
&&{}+2\sqrt{\epsilon}\int_0^t\<X^{\epsilon,h_{\epsilon}}(s)-X^{\epsilon,h_{\epsilon},\alpha}(s),(\sigma_{\epsilon}(X^{\epsilon,h_{\epsilon}}(s))-\sigma_{\epsilon}(X^{\epsilon,h_{\epsilon},\alpha}(s)))\dif W(s)\>_{\mR^m}\\
&&{}+\frac{\epsilon}{2}\int_0^t\|\sigma_{\epsilon}(X^{\epsilon,h_{\epsilon}}(s))-\sigma_{\epsilon}(X^{\epsilon,h_{\epsilon},\alpha}(s))\|_{L_2(l^2;\mR^m)}^2\dif s\\
&:=&\sum_{i=1}^7I_i^{\epsilon,\alpha}(t).
\de
By BDG's inequality, {\bf (H2)}, Lemma \ref{b1}-\ref{b2}, H\"{o}lder's inequality and Young's inequality, we have
\ce
&&\bE\sup_{s\in[0,t]}|I_5^{\epsilon,\alpha}(s)|\\
&\leq& C\bE[(\int_0^t|\widetilde{X}^{\epsilon,h_{\epsilon},\alpha}(s)-X^{\epsilon,h_{\epsilon}}(s)|^2\|\sigma_{\epsilon}(X^{\epsilon,h_{\epsilon}}(s))-\sigma_{\epsilon}(X^{\epsilon,h_{\epsilon},\alpha}(s))\|_{L_2(l^2;\mR^m)}^2\dif s)^{1/2}]\\
&\leq& C\bE[\sup_{s\in[0,t]}\|\sigma_{\epsilon}(X^{\epsilon,h_{\epsilon}}(s))-\sigma_{\epsilon}(X^{\epsilon,h_{\epsilon},\alpha}(s))\|_{L_2(l^2;\mR^m)}\sup_{s\in[0,t]}|\widetilde{X}^{\epsilon,h_{\epsilon},\alpha}(s)-X^{\epsilon,h_{\epsilon}}(s)|]\\
&\leq& C(\bE\sup_{s\in[0,t]}\|\sigma_{\epsilon}(X^{\epsilon,h_{\epsilon}}(s))-\sigma_{\epsilon}(X^{\epsilon,h_{\epsilon},\alpha}(s))\|_{L_2(l^2;\mR^m)}^2)^{1/2}(\bE\sup_{s\in[0,t]}|\widetilde{X}^{\epsilon,h_{\epsilon},\alpha}(s)-X^{\epsilon,h_{\epsilon}}(s)|^2)^{1/2}\\
&\leq& C(1+\bE\sup_{s\in[0,t]}|X^{\epsilon,h_{\epsilon}}(s)|+\bE\sup_{s\in[0,t]}|X^{\epsilon,h_{\epsilon},\alpha}(s)|)(\bE\sup_{s\in[0,t]}|\widetilde{X}^{\epsilon,h_{\epsilon},\alpha}(s)-X^{\epsilon,h_{\epsilon}}(s)|^2)^{1/2}\\
&\leq& C(\bE\sup_{s\in[0,t]}|\widetilde{X}^{\epsilon,h_{\epsilon},\alpha}(s)-X^{\epsilon,h_{\epsilon}}(s)|^2)^{1/2}
\de
and
\ce
&&\bE\sup_{s\in[0,t]}|I_6^{\epsilon,\alpha}(s)|\\
&\leq& C\bE[(\int_0^t|X^{\epsilon,h_{\epsilon}}(s)-X^{\epsilon,h_{\epsilon},\alpha}(s)|^2\|\sigma_{\epsilon}(X^{\epsilon,h_{\epsilon}}(s))-\sigma_{\epsilon}(X^{\epsilon,h_{\epsilon},\alpha}(s))\|_{L_2(l^2;\mR^m)}^2\dif s)^{1/2}]\\
&\leq&\frac{1}{4}\bE\sup_{s\in[0,t]}|X^{\epsilon,h_{\epsilon}}(s)-X^{\epsilon,h_{\epsilon},\alpha}(s)|^2\\
&&{}+C\int_0^t\bE[|X^{\epsilon,h_{\epsilon}}(s)-X^{\epsilon,h_{\epsilon},\alpha}(s)|^2(1\vee\log|X^{\epsilon,h_{\epsilon}}(s)-X^{\epsilon,h_{\epsilon},\alpha}(s)|^{-1})]\dif s.
\de
By H\"{o}lder's inequality, Young's inequality, {\bf (H2)}-{\bf (H3)} and Lemma \ref{b1}-\ref{b2}, we have
\ce
&&\bE\sup_{s\in[0,t]}|I_1^{\epsilon,\alpha}(s)|\\
&\leq& C(1+\bE\sup_{s\in[0,t]}|X^{\epsilon,h_{\epsilon}}(s)|+\bE\sup_{s\in[0,t]}|X^{\epsilon,h_{\epsilon},\alpha}(s)|)(\bE\sup_{s\in[0,t]}|\widetilde{X}^{\epsilon,h_{\epsilon},\alpha}(s)-X^{\epsilon,h_{\epsilon}}(s)|^2)^{1/2}\\
&\leq& C(\bE\sup_{s\in[0,t]}|\widetilde{X}^{\epsilon,h_{\epsilon},\alpha}(s)-X^{\epsilon,h_{\epsilon}}(s)|^2)^{1/2},
\de
\ce
\bE\sup_{s\in[0,t]}|I_2^{\epsilon,\alpha}(s)|\leq C\int_0^t\bE[|X^{\epsilon,h_{\epsilon}}(s)-X^{\epsilon,h_{\epsilon},\alpha}(s)|^2(1\vee\log|X^{\epsilon,h_{\epsilon}}(s)-X^{\epsilon,h_{\epsilon},\alpha}(s)|^{-1})]\dif s,
\de
\ce
&&\bE\sup_{s\in[0,t]}|I_3^{\epsilon,\alpha}(s)|\\
&\leq& C\bE[(\int_0^t|\widetilde{X}^{\epsilon,h_{\epsilon},\alpha}(s)-X^{\epsilon,h_{\epsilon}}(s)|^2\|\sigma_{\epsilon}(X^{\epsilon,h_{\epsilon}}(s))-\sigma_{\epsilon}(X^{\epsilon,h_{\epsilon},\alpha}(s))\|_{L_2(l^2;\mR^m)}^2\dif s)^{1/2}(\int_0^t\dot{h}_{\epsilon}(s)^2\dif s)^{1/2}]\\
&\leq& C(1+\bE\sup_{s\in[0,t]}|X^{\epsilon,h_{\epsilon}}(s)|+\bE\sup_{s\in[0,t]}|X^{\epsilon,h_{\epsilon},\alpha}(s)|)(\bE\sup_{s\in[0,t]}|\widetilde{X}^{\epsilon,h_{\epsilon},\alpha}(s)-X^{\epsilon,h_{\epsilon}}(s)|^2)^{1/2}\\
&\leq& C(\bE\sup_{s\in[0,t]}|\widetilde{X}^{\epsilon,h_{\epsilon},\alpha}(s)-X^{\epsilon,h_{\epsilon}}(s)|^2)^{1/2}
\de
and
\ce
&&\bE\sup_{s\in[0,t]}|I_4^{\epsilon,\alpha}(s)|\\
&\leq& C\bE[(\int_0^t|X^{\epsilon,h_{\epsilon}}(s)-X^{\epsilon,h_{\epsilon},\alpha}(s)|^2\|\sigma_{\epsilon}(X^{\epsilon,h_{\epsilon}}(s))-\sigma_{\epsilon}(X^{\epsilon,h_{\epsilon},\alpha}(s))\|_{L_2(l^2;\mR^m)}^2\dif s)^{1/2}(\int_0^t\dot{h}_{\epsilon}(s)^2\dif s)^{1/2}]\\
&\leq&\frac{1}{4}\bE\sup_{s\in[0,t]}|X^{\epsilon,h_{\epsilon}}(s)-X^{\epsilon,h_{\epsilon},\alpha}(s)|^2\\
&&{}+C\int_0^t\bE[|X^{\epsilon,h_{\epsilon}}(s)-X^{\epsilon,h_{\epsilon},\alpha}(s)|^2(1\vee\log|X^{\epsilon,h_{\epsilon}}(s)-X^{\epsilon,h_{\epsilon},\alpha}(s)|^{-1})]\dif s.
\de
We also have
\ce
\bE\sup_{s\in[0,t]}|I_7^{\epsilon,\alpha}(s)|\leq C\int_0^t\bE[|X^{\epsilon,h_{\epsilon}}(s)-X^{\epsilon,h_{\epsilon},\alpha}(s)|^2(1\vee\log|X^{\epsilon,h_{\epsilon}}(s)-X^{\epsilon,h_{\epsilon},\alpha}(s)|^{-1})]\dif s.
\de

Denote
\ce
f(t):=\bE\sup_{s\in[0,t]}|X^{\epsilon,h_{\epsilon},\alpha}(s)-X^{\epsilon,h_{\epsilon}}(s)|^2.
\de
Combining the above calculations and noting that there exists an $\eta>0$ such that
\ce
r^2(1\vee\log r^{-1})\leq\rho_{\eta}(r^2),
\de
we have
\be\label{bihari}
f(t)&\leq&2\bE\sup_{s\in[0,t]}|X^{\epsilon,h_{\epsilon},\alpha}(s)-\widetilde{X}^{\epsilon,h_{\epsilon},\alpha}(s)|^2+2\bE\sup_{s\in[0,t]}|\widetilde{X}^{\epsilon,h_{\epsilon},\alpha}(s)-X^{\epsilon,h_{\epsilon}}(s)|^2\nonumber\\
&\leq&\frac{1}{2}f(t)+C\int_0^t\rho_{\eta}(f(s))\dif s+C(\bE\sup_{s\in[0,t]}|\widetilde{X}^{\epsilon,h_{\epsilon},\alpha}(s)-X^{\epsilon,h_{\epsilon}}(s)|^2)^{1/2}\nonumber\\
&&{}+2\bE\sup_{s\in[0,t]}|\widetilde{X}^{\epsilon,h_{\epsilon},\alpha}(s)-X^{\epsilon,h_{\epsilon}}(s)|^2,
\ee
where we have used Jensen's inequality.
It is easy to see that $C$ in (\ref{bihari}) is independent of $\epsilon$ and $\alpha$. Therefore the Bihari's inequality and (\ref{key1}) yields
\ce
\lim_{\alpha\rightarrow0}\bE\sup_{t\in[0,T]}|X^{\epsilon,h_{\epsilon},\alpha}(t,x)-X^{\epsilon,h_{\epsilon}}(t,x)|^2=0
\de
uniformly in $\epsilon>0$, and we complete the proof.
\end{proof}

Similarly, let $(X^h,K^h)$ solve the following equation:
\ce
\begin{cases}
\dif X^h(t)\in b(X^h(t))\dif t+\sigma(X^h(t))\dot{h}(t)\dif t-A(X^h(t))\dif t,\\
X(0)=x\in \overline{D(A)},
\end{cases}
\de
and let us consider their related Yosida approximation
\ce
\begin{cases}
\dif X^{h,\alpha}(t)= b(X^{h,\alpha}(t))\dif t+\sigma(X^{h,\alpha}(t))\dot{h}(t)\dif t-A^{\alpha}(X^{h,\alpha}(t))\dif t,\\
X^{h,\alpha}(0)=x\in \overline{D(A)}.
\end{cases}
\de

Analogously, we can prove the following lemma.
\bl\label{p2}
Assume that {\bf (H1)}-{\bf (H4)} hold. For any $p\geq 1$ and $x\in \overline{D(A)}$, there exists $C>0$ such that for any $\alpha>0$,
\ce
\sup_{t\in[0,T]}|X^{h,\alpha}(t,x)|^{2p}\leq C,
\de
where $C$ may be dependent of $p$, $T$, $N$ and $x$, but is independent of $\epsilon$ and $\alpha$.

Moreover, for any $x\in \overline{D(A)}$, we have
\ce
\lim_{\alpha\rightarrow0}\sup_{t\in[0,T]}|X^{h,\alpha}(t,x)-X^h(t,x)|^2=0
\de
uniformly in $\epsilon>0$.
\el

\subsection{Main estimates}

We first give the following uniform estimate.
\bl\label{uniform estimate}
Assume that {\bf (H1)}-{\bf (H3)} hold. Then, for any $p\geq 1$, there exists $C>0$ such that for any $\epsilon\in(0,1)$ and $x,y\in \overline{D(A)}$
\ce
\bE\sup_{t\in[0,T]}|X^{\epsilon,h_{\epsilon}}(t,x)-X^{\epsilon,h_{\epsilon}}(t,y)|^{2p}\leq C|x-y|^{2p},
\de
where $C$ may be dependent of $p$, $T$ and $N$, but is independent of $\epsilon$.
\el
\begin{proof}
Noting that there exists an $\eta>0$ such that
\ce
r^2(1\vee\log r^{-1})\leq\rho_{\eta}(r^2),
\de
and applying Bihari's inequality, this result actually follows similarly as in the proof of \cite[Lemma 3.4]{RXZ}.
\end{proof}

The following lemma states a square convergence result.
\bl\label{p3}
Assume that {\bf (H2)}-{\bf (H4)} hold. Then for every $\alpha>0$ and every $x\in \overline{D(A)}$, we have
\ce
\lim_{\epsilon\rightarrow0}\bE\sup_{t\in[0,T]}|X^{\epsilon,h_{\epsilon},\alpha}(t,x)-X^{h,\alpha}(t,x)|^2=0.
\de
\el
\begin{proof}
Denote
\ce
u_{\epsilon}(t)&:=&X^{\epsilon,h_{\epsilon},\alpha}(t,x)-X^{h,\alpha}(t,x)\\
&=&\int_0^t(b_{\epsilon}(X^{\epsilon,h_{\epsilon},\alpha}(s))-b(X^{h,\alpha}(s)))\dif s+\int_0^t(\sigma_{\epsilon}(X^{\epsilon,h_{\epsilon},\alpha}(s))\dot{h}_{\epsilon}-\sigma(X^{h,\alpha}(s))\dot{h}(s))\dif s\\
&&{}+\sqrt{\epsilon}\int_0^t\sigma_{\epsilon}(X^{\epsilon,h_{\epsilon},\alpha}(s))\dif W(s)+\int_0^t(A^{\alpha}(X^{h,\alpha}(s))-A_{\epsilon}^{\alpha}(X^{\epsilon,h_{\epsilon},\alpha}(s)))\dif s.
\de
By It\^{o}'s formula, we have
\ce
|u_{\epsilon}(t)|^2
&=&2\int_0^t\<u_{\epsilon}(s),b_{\epsilon}(X^{\epsilon,h_{\epsilon},\alpha}(s))-b(X^{\epsilon,h_{\epsilon},\alpha}(s))\>_{\mR^m}\dif s\\
&&{}+2\int_0^t\<u_{\epsilon}(s),b(X^{\epsilon,h_{\epsilon},\alpha}(s))-b(X^{h,\alpha}(s))\>_{\mR^m}\dif s\\
&&{}+2\int_0^t\<u_{\epsilon}(s),\sigma_{\epsilon}(X^{\epsilon,h_{\epsilon},\alpha}(s))\dot{h}_{\epsilon}(s)-\sigma(X^{\epsilon,h_{\epsilon},\alpha}(s))\dot{h}_{\epsilon}(s)\>_{\mR^m}\dif s\\
&&{}+2\int_0^t\<u_{\epsilon}(s),(\sigma(X^{\epsilon,h_{\epsilon},\alpha}(s))-\sigma(X^{h,\alpha}(s)))\dot{h}_{\epsilon}(s)\>_{\mR^m}\dif s\\
&&{}+2\int_0^t\<u_{\epsilon}(s),\sigma(X^{h,\alpha}(s))(\dot{h}_{\epsilon}(s)-\dot{h}(s))\>_{\mR^m}\dif s\\
&&{}+2\sqrt{\epsilon}\int_0^t\<u_{\epsilon}(s),\sigma_\epsilon(X^{\epsilon,h_{\epsilon},\alpha}(s))\dif W(s)\>_{\mR^m}+\epsilon\int_0^t\|\sigma_{\epsilon}(X^{\epsilon,h_{\epsilon},\alpha}(s))\|_{L_2(l^2;\mR^m)}^2\dif s\\
&&{}+2\int_0^t\<u_{\epsilon}(s),A^{\alpha}(X^{h,\alpha}(s))-A^{\alpha}(X^{\epsilon,h_{\epsilon},\alpha}(s))\>_{\mR^m}\\
&&{}+2\int_0^t\<u_{\epsilon}(s),A^{\alpha}(X^{\epsilon,h_{\epsilon},\alpha}(s))-A_{\epsilon}^{\alpha}(X^{\epsilon,h_{\epsilon},\alpha}(s))\>_{\mR^m}\\
&:=&\sum_{i=1}^9I_i^{\epsilon,\alpha}(t).
\de
By Young's inequality, we have
\ce
\bE\sup_{t\in[0,T]}|I_1^{\epsilon,\alpha}(t)|\leq \frac{1}{8}\bE\sup_{t\in[0,T]}|u_{\epsilon}(t)|^2+C\int_0^T\bE[|b_{\epsilon}(X^{\epsilon,h_{\epsilon},\alpha}(t))-b(X^{\epsilon,h_{\epsilon},\alpha}(t))|^2]\dif t,
\de
\ce
\bE\sup_{t\in[0,T]}|I_3^{\epsilon,\alpha}(t)|\leq \frac{1}{8}\bE\sup_{t\in[0,T]}|u_{\epsilon}(t)|^2+C\int_0^T\bE[\|\sigma_{\epsilon}(X^{\epsilon,h_{\epsilon},\alpha}(t))-\sigma(X^{\epsilon,h_{\epsilon},\alpha}(t))\|_{L_2(l^2;\mR^m)}^2]\dif t,
\de
\ce
\bE\sup_{t\in[0,T]}|I_9^{\epsilon,\alpha}(t)|\leq \frac{1}{8}\bE\sup_{t\in[0,T]}|u_{\epsilon}(t)|^2+C\int_0^T\bE[|A^{\alpha}(X^{\epsilon,h_{\epsilon},\alpha}(t))-A_{\epsilon}^{\alpha}(X^{\epsilon,h_{\epsilon},\alpha}(t))|^2]\dif t.
\de
Similarly, using {\bf (H2)} and {\bf (H3)}, we also have
\ce
\bE\sup_{t\in[0,T]}|I_2^{\epsilon,\alpha}(t)|\leq C\int_0^T\bE[|u_{\epsilon}(t)|^2(1\vee\log|u_{\epsilon}(t)|^{-1})]\dif t,
\de
\ce
\bE\sup_{t\in[0,T]}|I_4^{\epsilon,\alpha}(t)|\leq \frac{1}{8}\bE\sup_{t\in[0,T]}|u_{\epsilon}(t)|^2+C\int_0^T\bE[|u_{\epsilon}(t)|^2(1\vee\log|u_{\epsilon}(t)|^{-1})]\dif t,
\de
and by the monotonicity of $A^\alpha$,
\ce
I_8^{\epsilon,\alpha}(t)\leq0.
\de
By BDG's inequality, {\bf (H2)}-{\bf (H3)}, Lemma \ref{b1} and Lemma \ref{p2}, we have
\ce
\bE\sup_{t\in[0,T]}|I_6^{\epsilon,\alpha}(t)|+\bE\sup_{t\in[0,T]}|I_7^{\epsilon,\alpha}(t)|\leq C\sqrt{\epsilon}.
\de
Consequently, combining the above estimates and noting that there exists an $\eta>0$ such that
\ce
r^2(1\vee\log r^{-1})\leq\rho_{\eta}(r^2),
\de
we obtain
\ce
\bE\sup_{t\in[0,T]}|u_{\epsilon}(s)|^2&\leq& \frac{1}{2}\bE\sup_{t\in[0,T]}|u_{\epsilon}(t)|^2+C\int_0^T\bE[\rho_{\eta}(\sup_{s\in[0,t]}|u_{\epsilon}(s)|^2)]\dif t+C\sqrt{\epsilon}\\
&&{}+\bE\sup_{t\in[0,T]}|I_5^{\epsilon,\alpha}(t)|
+\int_0^T\bE[|b_{\epsilon}(X^{\epsilon,h_{\epsilon},\alpha}(t))-b(X^{\epsilon,h_{\epsilon},\alpha}(t))|^2]\dif t\\
&&{}+\int_0^T\bE[\|\sigma_{\epsilon}(X^{\epsilon,h_{\epsilon},\alpha}(t))-\sigma(X^{\epsilon,h_{\epsilon},\alpha}(t))\|_{L_2(l^2;\mR^m)}^2]\dif t\\
&&{}+\int_0^T\bE[|A^{\alpha}(X^{\epsilon,h_{\epsilon},\alpha}(t))-A_{\epsilon}^{\alpha}(X^{\epsilon,h_{\epsilon},\alpha}(t))|^2]\dif t\\
&\leq&\frac{1}{2}\bE\sup_{t\in[0,T]}|u_{\epsilon}(t)|^2+C\int_0^T\rho_{\eta}(\bE\sup_{s\in[0,t]}|u_{\epsilon}(s)|^2)\dif t+C\sqrt{\epsilon}\\
&&{}+\bE\sup_{t\in[0,T]}|I_5^{\epsilon,\alpha}(t)|+\int_0^T\bE[|b_{\epsilon}(X^{\epsilon,h_{\epsilon},\alpha}(t))-b(X^{\epsilon,h_{\epsilon},\alpha}(t))|^2]\dif t\\
&&{}+\int_0^T\bE[\|\sigma_{\epsilon}(X^{\epsilon,h_{\epsilon},\alpha}(t))-\sigma(X^{\epsilon,h_{\epsilon},\alpha}(t))\|_{L_2(l^2;\mR^m)}^2]\dif t\\
&&{}+\int_0^T\bE[|A^{\alpha}(X^{\epsilon,h_{\epsilon},\alpha}(t))-A_{\epsilon}^{\alpha}(X^{\epsilon,h_{\epsilon},\alpha}(t))|^2]\dif t,
\de
where we have used Jensen's inequality.
By an arguments similar to the proof of \cite[Lemma 3.7]{RXZ}, together with Lemma \ref{b1} and Lemma \ref{p2}, we can prove that
\ce
\lim_{\epsilon\rightarrow0}\bE\sup_{t\in[0,T]}|I_5^{\epsilon,\alpha}(t)|=0.
\de
Next we show that
\be\label{y1}
\eta_\epsilon:=\int_0^T\bE[|A^{\alpha}(X^{\epsilon,h_{\epsilon},\alpha}(t))-A_{\epsilon}^{\alpha}(X^{\epsilon,h_{\epsilon},\alpha}(t))|^2]\dif t\to 0\quad\text{as}\quad\epsilon\downarrow0.
\ee
In fact, for all $k$ we have
\ce
\eta_\epsilon&=&
\int_0^T\bE[|A^{\alpha}(X^{\epsilon,h_{\epsilon},\alpha}(t))-
A_{\epsilon}^{\alpha}(X^{\epsilon,h_{\epsilon},\alpha}(t))|^21_{\{\sup_{t\in [0,T]}|X^{\epsilon,h_{\epsilon},\alpha}(t)|\leq k\}}]
\dif t\\
&&{}+\int_0^T\bE[|A^{\alpha}(X^{\epsilon,h_{\epsilon},\alpha}(t))-
A_{\epsilon}^{\alpha}(X^{\epsilon,h_{\epsilon},\alpha}(t))|^21_{\{\sup_{t\in [0,T]}|X^{\epsilon,h_{\epsilon},\alpha}(t)|> k\}}]
\dif t\\
&=:& I_1+I_2.
\de
Note that by the Lipschitz continuity of Yosida approximation and {\bf (H1)} we have
\ce
|A^{\alpha}(X^{\epsilon,h_{\epsilon},\alpha}(t))|+
|A_{\epsilon}^{\alpha}(X^{\epsilon,h_{\epsilon},\alpha}(t))|\leq C_\alpha(1+|X^{\epsilon,h_{\epsilon},\alpha}(t)|).
\de
Hence for $k\geq 1$.
\ce
I_2&\leq&\int_0^T\mathbf{E}[|A^{\alpha}(X^{\epsilon,h_{\epsilon},\alpha}(t))-
A_{\epsilon}^{\alpha}(X^{\epsilon,h_{\epsilon},\alpha}(t))|^21_{\{\sup_{t\in [0,T]}|X^{\epsilon,h_{\epsilon},\alpha}(t)|> k\}}]dt\\
&\leq&C\mathbf{E}[\sup_{t\in[0,T]}|X^{\epsilon,h_{\epsilon},\alpha}(t)|^21_{\{\sup_{t\in [0,T]}|X^{\epsilon,h_{\epsilon},\alpha}(t)|> k\}}]\\
&\leq&C\mathbf{E}[\sup_{t\in[0,T]}|X^{\epsilon,h_{\epsilon},\alpha}(t)|^2\frac{\sup_{t\in [0,T]}|X^{\epsilon,h_{\epsilon},\alpha}(t)|}{k}]\\
&\leq&C\sup_\epsilon\frac{\mathbf{E}[\sup_{t\in [0,T]}|X^{\epsilon,h_{\epsilon},\alpha}(t)|^3]}{k}.
\de
Thus by letting first $\epsilon\downarrow 0$ and then $k\to \infty$, we have (\ref{y1}).
Similarly,
\ce
&&\lim_{\epsilon\rightarrow0}\int_0^T\bE[|b_{\epsilon}(X^{\epsilon,h_{\epsilon},\alpha}(t))-b(X^{\epsilon,h_{\epsilon},\alpha}(t))|^2]\dif t=0,\\
&&\lim_{\epsilon\rightarrow0}\int_0^T\bE[\|\sigma_{\epsilon}(X^{\epsilon,h_{\epsilon},\alpha}(t))-\sigma(X^{\epsilon,h_{\epsilon},\alpha}(t))\|_{L_2(l^2;\mR^m)}^2]\dif t=0.
\de
Now the proof is completed by Bihari's inequality.
\end{proof}

\subsection{Completion of the proof of Theorem \ref{Main}}

Since
\ce
\bE\sup_{t\in[0,T]}|X^{\epsilon,h_{\epsilon}}(t,x)-X^h(t,x)|^2&\leq&3\bE\sup_{t\in[0,T]}|X^{\epsilon,h_{\epsilon}}(t,x)-X^{\epsilon,h_{\epsilon},\alpha}(t,x)|^2\\
&&{}+3\bE\sup_{t\in[0,T]}|X^{\epsilon,h_{\epsilon},\alpha}(t,x)-X^{h,\alpha}(t,x)|^2\\
&&{}+3\bE\sup_{t\in[0,T]}|X^{h,\alpha}(t,x)-X^h(t,x)|^2,
\de
by Lemmas \ref{p1}, \ref{p2} and \ref{p3}, given $\xi>0$ we first choose $\alpha$, independent of $\epsilon$, such that the first and the third terms are less then $\xi/9$. Then having fixed $\alpha$ this way, we can choose $\epsilon$ such that the second term is less than $\xi/9$, and finally we have
\ce
\mE[\sup_{t\in[0,T]}|X^{\epsilon,h_{\epsilon}}(t,x)-X^h(t,x)|^2]\leq \xi.
\de
This implies for all $x\in \overline{D(A)}$,
\ce
\sup_{t\in[0,T]}|X^{\epsilon,h_{\epsilon}}(t,x)-X^h(t,x)|\rightarrow 0,\quad \text{in probability}.
\de
As in the proof of \cite[Lemma 3.7]{RXZ}, we can strengthen it by Lemma \ref{uniform estimate}
\be\label{key estimate}
\xi_{n,\epsilon}:=\sup_{t\in[0,T],x\in\overline{D(A)},|x|\leq n}|X^{\epsilon,h_{\epsilon}}(t,x)-X^h(t,x)|\rightarrow 0,\quad\text{in probability}.
\ee

As before, $h_\eps\in\cA^T_N$ converge almost surely to $h\in\cA^T_N$
as random variables in $\ell^2_T$. Since $\cD_N$ is compact and the law of $W$ is tight, $\{h_\eps,W\}$ is tight in $\cD_N\times\mathcal {C}_T(\mU)$. We assume that the law of $\{h_\eps,W\}$ weakly converges to $\mu$. Then the law of $h$ is just $\mu(\cdot,\mathcal {C}_T(\mU))$. By Skorohod's representation theorem, there exist another probability space $(\tilde \Omega,\tilde P)$ and
$\{\tilde h_\eps,\tilde W^\eps\}$ and $\{\tilde h,\tilde W\}$ on it such that

(1) $(\tilde h_\eps,\tilde W^\eps)$ a.s. converges to $(\tilde h,\tilde W)$;

(2) $(\tilde h_\eps,\tilde W^\eps)$ has the same law as $(h_\eps,W)$;

(3) The law of $\{\tilde h,\tilde W\}$ is $\mu$, and the law of
$h$ is the same as $\tilde h$.

By the classical Yamada-Watanabe theorem, for the solution $X^{\epsilon,h_{\epsilon}}$ of Eq. (\ref{control}), we have
\ce
X^{\epsilon,h_{\epsilon}}=\Phi(\frac{1}{\sqrt{\eps}}\int^\cdot_0\dot{h}_\eps\dif s+W^\eps),
\de
where $\Phi$ is the strong solution functional (cf. \cite{IW}). Therefore applying (\ref{key estimate}), we get
\ce
\Phi(\frac{1}{\sqrt{\eps}}\int^\cdot_0\dot{\tilde h}_\eps\dif s+\tilde W^\eps)\to X^{\tilde h}, \quad\text{in probability},
\de
from which we can immediately have
\ce
\Phi(\frac{1}{\sqrt{\eps}}\int^\cdot_0\dot{h}_\eps\dif s+W)\to X^{h}, \quad\text{in distribution}.
\de
Thus, {\bf (LD)$_\mathbf{1}$} holds.

{\bf (LD)$_\mathbf{2}$} can be verified in a similar way.

Thus, by Theorem \ref{Th2}, we have proved Theorem \ref{Main}.

\subsection{An Example.}

Consider the following stochastic variational inequalities
\be\label{svi1}
\left\{
\begin{array}{ll}
\dif X^\eps(t)\in b_\eps(X^\eps(t))\dif t+\sqrt{\eps}
\sigma_\eps(X^\eps(t))\dif W(t)-\partial\vph_\eps(X^\eps(t))\dif t,\\
X^\eps(0)=x\in \overline{D(\varphi)},\ \ \eps\in(0,1],
\end{array}
\right.
\ee
where $\vph_\eps$ and $\vph$ are lower semicontinuous convex functions with common domain $D$ whose interior is non-vide. Suppose $\lim_{\eps\downarrow 0}\vph_\eps(x)=\vph(x)$ for every $x\in D$. Define the Moreau-Yosida approximation of $\vph$ by
\ce
\vph^\alpha(x)=\inf_{y\in \mR^m}\{\vph(y)+|y-x|^2/{(2\alpha)}\}
\de
and that of $\vph_\eps$ in the same way. Then, by \cite[Ch. II]{Ba}, $\vph^\alpha_\eps$ and $\vph^\alpha$ are differentiable convex functions defined on the whole space and it is not difficult to prove that
\ce
\lim_{\eps\downarrow 0}\vph_\eps^\alpha(x)=\vph^\alpha(x)
\de
for every $\alpha$. Moreover, denoting
\ce
A:=\partial \vph,\quad A_\eps:=\partial \vph_\eps,
\de
we have
\ce
A^\alpha=\partial \vph^\alpha,\quad A^\alpha_\eps=\partial \vph_\eps^\alpha.
\de
Hence by \cite[Ch.D, Corollary 6.2.7 and Corollary 6.2.8]{HL} {\bf (H1)} and {\bf (H4)} are satisfied by $A$ and $A_\eps$. Consequently, if $\sigma$ and $b$ satisfy {\bf (H2)} and {\bf (H3)} respectively, then Theorem \ref{Main} is applicable.

\section{Functional iterated logarithm law}

In this section, we derive a functional iterated logarithm law.

Let now $W$ be a $k$-dimensional Brownian motion and $T=1$ be fixed. Let $U$ be an open subset of $\mR^m$ and we denote by $\mathcal {C}([0,1];U)$ the set of continuous path $u:[0,1]\rightarrow U$, endowed with the norm
\ce
\|u\|=\sup_{t\in[0,1]}|u(t)|.
\de
For $u,v\in \mathcal {C}([0,1];U)$, we denote
\ce
&&\|u\|=\sup_{t\in[0,1]}|u(t)|,\quad d(u,v)=\|u-v\|;\\
&&d(u,B)=\inf_{v\in B}d(u,v),\quad B\in \cB(\mathcal {C}([0,1];U)).
\de

Let
\ce
\cH_m^2:=\{u\in \mathcal {C}([0,1];\mR^m); u(0)=0,~\int_0^1|\dot{u}(s)|^2\dif s<\infty\}.
\de
 Then $\cH_m^2$ is a Hilbert space endowed with the following scalar product
\ce
\<u,v\>:=\int_0^1\<\dot{u}(t),\dot{v}(s)\>_{\mR^m}\dif s.
\de
We set
\ce
\|u\|_{\cH_m^2}^2=\<u,u\>^{1/2}=\|\dot{u}\|_{L^2([0,1];\mR^m)}.
\de

We first state the following continuous dependence result which will be needed in the proof of Theorem \ref{exponential}.

\bl\label{continuity}
Suppose that $\widetilde{A}$ is a maximal monotone operator with $0\in\mathrm{Int}(D(\widetilde{A}))$. Suppose further that $\tilde{\sigma}$ and $\tilde{b}$ are continuous functions and satisfy that for some $C>0$ and all $x,y\in\mR^m$
\be\label{assumption1}
\|\tilde{\sigma}(x)-\tilde{\sigma}(y)\|_{L_2(l^2;\mR^m)}^2\vee
\<x-y,\tilde{b}(x)-\tilde{b}(y)\>\leq C|x-y|^2(1\vee\log|x-y|^{-1})
\ee
and
\be\label{assumption2}
\|\tilde{\sigma}(x)\|_{L_2(l^2;\mR^m)}\vee|\tilde{b}(x)|\leq C(1+|x|).
\ee
Then $B_R\ni h\rightarrow Z^h$ is continuous with respect to the distance $d$,
where $B_R:=\{h\in\cH_k^2;\frac{1}{2}\|h\|_{\cH_k^2}^2\leq 2R\}$ and $(Z^h(\cdot,x),K^h(\cdot,x))$ solves the following MSDE:
\ce
\begin{cases}
\dif Z^h(t)\in\tilde{b}(Z^h(t))\dif t+\tilde{\sigma}(Z^h(t))\dot{h}(t)\dif t-\widetilde{A}(Z^h(t))\dif t,\\
Z^h(0)=x\in\overline{D(\widetilde{A})}.
\end{cases}
\de
\el
\begin{proof}
For fixed $h_0\in B_R$, let $h\in B_R$ such that
\ce
\sup_{t\in[0,1]}|h(t)-h_0(t)|\rightarrow0.
\de
Since
\ce
|Z^h(t)-Z^{h_0}(t)|^2&=&2\int_0^t\<Z^h(s)-Z^{h_0}(s),\tilde{b}(Z^h(s))-\tilde{b}(Z^{h_0}(s))\>_{\mR^m}\dif s\\
&&+2\int_0^t\<Z^h(s)-Z^{h_0}(s),\tilde{\sigma}(Z^h(s))\dot{h}(s)-\tilde{\sigma}(Z^{h_0}(s))\dot{h}_0(s)\>_{\mR^m}\dif s\\
&&-2\int_0^t\<Z^h(s)-Z^{h_0}(s),\dif K^h(s)-\dif K^{h_0}(s)\>_{\mR^m}\\
&\leq&2\int_0^t\<Z^h(s)-Z^{h_0}(s),\tilde{b}(Z^h(s))-\tilde{b}(Z^{h_0}(s))\>_{\mR^m}\dif s\\
&&+2\int_0^t\<Z^h(s)-Z^{h_0}(s),\tilde{\sigma}(Z^h(s))\dot{h}(s)-\tilde{\sigma}(Z^{h_{0}}(s))\dot{h}(s)\>_{\mR^m}\dif s\\
&&+2\int_0^t\<Z^h(s)-Z^{h_0}(s),\tilde{\sigma}(Z^{h_0}(s))\dot{h}(s)-\tilde{\sigma}(Z^{h_0}(s))\dot{h}_0(s)\>_{\mR^m}\dif s,
\de
by Young's inequality, (\ref{assumption1}) and the fact that there exists an $\eta>0$ satisfying
\ce
r^2(1\vee\log r^{-1})\leq\rho_{\eta}(r^2),
\de
we have
\ce
\sup_{s\in[0,t]}|Z^h(s)-Z^{h_0}(s)|^2&\leq&2\int_0^t|Z^h(s)-Z^{h_0}(s)|^2(1\vee\log|Z^h(s)-Z^{h_0}(s)|^{-1})\dif s\\
&&+2\int_0^t|Z^h(s)-Z^{h_0}(s)|\cdot\|\tilde{\sigma}(Z^h(s))-\tilde{\sigma}(Z^{h_0}(s))\|_{L_2(l^2;\mR^m)}\cdot|\dot{h}(s)|\dif s\\
&&+\sup_{t'\in[0,t]}2\int_0^{t'}\<Z^h(s)-Z^{h_0}(s),\tilde{\sigma}(Z^{h_0}(s))\dot{h}(s)-\tilde{\sigma}(Z^{h_0}(s))\dot{h}_0(s)\>_{\mR^m}\dif s\\
&\leq&2\int_0^t|Z^h(s)-Z^{h_0}(s)|^2(1\vee\log|Z^h(s)-Z^{h_0}(s)|^{-1})\dif s\\
&&+4R^{1/2}(\int_0^t|Z^h(s)-Z^{h_0}(s)|^2\cdot\|\tilde{\sigma}(Z^h(s))-\tilde{\sigma}(Z^{h_0}(s))\|_{L_2(l^2;\mR^m)}^2\dif s)^{1/2}\\
&&+\sup_{t'\in[0,t]}2\int_0^{t'}\<Z^h(s)-Z^{h_0}(s),\tilde{\sigma}(Z^{h_0}(s))(\dot{h}(s)-\dot{h}_0(s))\>_{\mR^m}\dif s\\
&\leq&\frac{1}{2}\sup_{s\in[0,t]}|Z^h(s)-Z^{h_0}(s)|^2+C\int_0^t\rho_{\eta}(\sup_{t'\in[0,s]}|Z^h(t')-Z^{h_0}(t')|^2)\dif s\\
&&+\sup_{t'\in[0,t]}2\int_0^{t'}\<Z^h(s)-Z^{h_0}(s),\tilde{\sigma}(Z^{h_0}(s))(\dot{h}(s)-\dot{h}_0(s))\>_{\mR^m}\dif s.
\de

Denote
\ce
I:=\sup_{t'\in[0,t]}2\int_0^{t'}\<Z^h(s)-Z^{h_0}(s),\tilde{\sigma}(Z^{h_0}(s))(\dot{h}(s)-\dot{h}_0(s))\>_{\mR^m}\dif s.
\de
Define
\ce
u(t):=\int_0^t\tilde{\sigma}(Z^{h_0}(s))(\dot{h}(s)-\dot{h}_0(s))\dif s,\quad
v(t):=Z^h(t)-Z^{h_0}(t).
\de
Integration by parts gives
\ce
&&\int_0^t\<Z^h(s)-Z^{h_0}(s),\tilde{\sigma}(Z^{h_0}(s))(\dot{h}(s)-\dot{h}_0(s))\>_{\mR^m}\dif s\\
&=&\<u(t),v(t)\>_{\mR^m}-\int_0^t\<u(s),\dif K^h(s)-\dif K^{h_0}(s)\>_{\mR^m}\\
&&{}-\int_0^t\<u(s),\tilde{b}(Z^h(s))-\tilde{b}(Z^{h_0}(s))\>_{\mR^m}\dif s\\
&&{}-\int_0^t\<u(s),\tilde{\sigma}(Z^h(s))\dot{h}(s)-\tilde{\sigma}(Z^{h_0}(s))\dot{h}_0(s)\>_{\mR^m}\dif s\\
&:=&I_1(t)+I_2(t)+I_3(t)+I_4(t).
\de
By the proof of Lemma \ref{b2}, we have for $h,h_0\in B_R$
\ce
&&\sup_{t\in[0,1]}|Z^h(t)|^p+|K^h(\cdot)|_1^0\leq C_R,\\
&&
\sup_{t\in[0,1]}|Z^{h_0}(t)|^p+|K^{h_0}(\cdot)|_1^0\leq C_R.
\de
Thus by assumptions (\ref{assumption1}) and (\ref{assumption2}),
\ce
\sup_{t\in[0,1]}(|I_1(t)|+|I_2(t)|+|I_3(t)|+|I_4(t)|)\leq C\sup_{t\in[0,1]}|u(t)|.\\
\de
If we can show that
\be\label{key step}
\sup_{t\in[0,1]}|u(t)|\rightarrow0\quad\text{as}\quad
\sup_{t\in[0,1]}|h(t)-h_0(t)|\rightarrow0,
\ee
then we can deduce that $I\rightarrow0$ as $\sup_{t\in[0,1]}|h(t)-h_0(t)|\rightarrow0$ which in turn implies by Bihari's inequality that
\ce
\sup_{t\in[0,1]}|Z^h(t)-Z^{h_0}(t)|\rightarrow0\quad\text{as}\quad
\sup_{t\in[0,1]}|h(t)-h_0(t)|\rightarrow0,
\de
as desired.

Now we show (\ref{key step}).
Let $\rho$ be a mollifier with $supp~\rho\subset(-1,1)$, $\rho\in \mathcal {C}^{\infty}$ and $\int_{-1}^1\rho(t)\dif t=1$. We denote $\rho_n(u):=n\rho(nu)$. Consider the following function:
\ce
\phi_n(t)=\int_0^1\tilde{\sigma}(Z^{h_0}(u))\rho_n(t-u)\dif u.
\de
It is obvious that $s\rightarrow\phi_n(s)$ is continuously differentiable and
\ce
\int_0^1\|\tilde{\sigma}(Z^{h_0}(t))-\phi_n(t)\|^2\dif t\rightarrow0\quad\text{as}\quad n\rightarrow\infty.
\de
Therefore
\ce
\sup_{t\in[0,1]}|u(t)|&=&\sup_{t\in[0,1]}|\int_0^t\tilde{\sigma}(Z^{h_0}(s))(\dot{h}(s)-\dot{h}_0(s))\dif s|\\
&\leq&\sup_{t\in[0,1]}|\int_0^t(\tilde{\sigma}(Z^{h_0}(s))-\phi_n(s))(\dot{h}(s)-\dot{h}_0(s))\dif s|\\
&&{}+\sup_{t\in[0,1]}|\int_0^t\phi_n(s)(\dot{h}(s)-\dot{h}_0(s))\dif s|=:J_1+J_2.
\de
By H\"{o}lder's inequality, for any $\delta>0$, there exists $n_0\in\mN$ such that
\ce
J_1\leq(\int_0^1\|\tilde{\sigma}(Z^h(t))-\phi_{n_0}(t)\|^2\dif t)^{1/2}(\int_0^1|\dot{h}(t)-\dot{h}_0(t)|^2\dif t)^{1/2}\leq4R^{1/2}\delta^{1/2}.
\de
For the second term $J_2$, integrating by parts we have
\ce
J_2&\leq&\sup_{t\in[0,1]}|h(t)-h_0(t)|(\sup_{t\in[0,1]}\|\phi_{n_0}(t)\|+\int_0^1\|\dot{\phi}_{n_0}(t)\|\dif t)\\
&\rightarrow&0,\quad \mbox{as}\quad \sup_{t\in[0,1]}|h(t)-h_0(t)|\rightarrow0,
\de
which is what we wanted and the proof is completed.
\end{proof}

The following result establishes a Freidlin-Wentzell estimate for $X^\eps$. In the case of ordinary stochastic differential equations, it is proved in \cite{BS}.

\bt\label{exponential}
Assume that {\bf (H1)}-{\bf (H4)} hold. Let $a>0$ be fixed. Then for every $\alpha>0$ and $R>0$, there exist $\epsilon_0,~\eta_0>0$ such that for every $h$ with $I(f)\leq a$ and $f=X^h$, we have
\ce
\bP\{d(\epsilon W,h)<\eta,d(X^{\epsilon},f)>\alpha\}\leq\exp(\frac{-R}{\epsilon^2})
\de
for every $0<\epsilon<\epsilon_0$ and $0<\eta<\eta_0$.
\et
\begin{proof}
Having Lemma \ref{continuity} in hand, the proof is similar to \cite{BS}, and we include it here for reader's convenience.

Consider the diffusion
\ce
Z^{\epsilon}(t,x)={\epsilon W(t)\choose X^{\epsilon}(t,x)},
\de
where $Z^{\epsilon}(t,x)$ is the solution of the following MSDE:
\ce
\begin{cases}
\dif Z^{\epsilon}(t)\in\tilde{b}_{\epsilon}(Z^{\epsilon}(t))\dif t+\epsilon\tilde{\sigma}_{\epsilon}(Z^{\epsilon}(t))\dif W_t-\widetilde{A}_{\epsilon}(Z^{\epsilon}(t))\dif t,\\
Z^{\epsilon}(0)={0\choose x},
\end{cases}
\de
where
\ce
\tilde{b}_{\epsilon}(x)={0\choose b_{\epsilon}(x)}:\mR^{k+m}\rightarrow\mR^{k+m},\quad
\tilde{\sigma}_{\epsilon}(x)={I\choose\sigma_{\epsilon}(x)}:\mR^{k+m}\rightarrow\mR^{(k+m)\times k}
\de
and
\ce
\widetilde{A}_{\epsilon}(x)={0\choose A_{\epsilon}(x)}:\mR^{k+m}\rightarrow2^{\mR^{k+m}}.
\de
For each $f\in C([0,1];\mR^{k+m})$, we define
\ce
\widetilde{I}(f)=\frac{1}{2}\inf_{\{h\in\cH_k^2;f=Z^h\}}\|h\|^2_{\cH_k^2},
\de
where $(Z^h(\cdot,x),K^h(\cdot,x))$ solves the following equation:
\be\label{compose}
\begin{cases}
\dif Z^h(t)\in\tilde{b}(Z^h(t))\dif t+\tilde{\sigma}(Z^h(t))\dot{h}(t)\dif t-\widetilde{A}(Z^h(t))\dif t,\\
Z^h(0)={0\choose x},
\end{cases}
\ee
where
\ce
\tilde{b}(x)={0\choose b(x)}:\mR^{k+m}\rightarrow\mR^{k+m},
\quad
\tilde{\sigma}(x)={I\choose\sigma(x)}:\mR^{k+m}\rightarrow\mR^{(k+m)\times k}
\de
and
\ce
\widetilde{A}(x)={0\choose A(x)}:\mR^{k+m}\rightarrow2^{\mR^{k+m}}.
\de
Obvious that $\tilde{b}_{\epsilon}$, $\tilde{b}$, $\tilde{\sigma}_{\epsilon}$, $\tilde{\sigma}$, $\widetilde{A}_{\epsilon}$ and $\widetilde{A}$ also satisfy assumptions {\bf (H1)}-{\bf (H4)}.
One can easily see that Eq. (\ref{compose}) is decomposed into two systems, i.e.,
\ce
Z^h={h\choose X^h}
\de
where $X^h(t,x)$ solves the following MDE:
\ce
\dif X^h(t)\in b(X^h(t))\dif t+
\sigma(X^h(t))\dot h(t)\dif t-A(X^h(t))\dif t,\quad
X^h(0)=x.
\de

Fix now $\alpha>0$ and $R>0$. Then
\ce
\bP\{d(X^{\epsilon},X^h)>\alpha,d(\epsilon W,h)<\eta\}=\bP\{Z^{\epsilon}\in\Lambda_{\eta}\},
\de
where
\ce
\Lambda_{\eta}=\{\chi={\chi_1\choose\chi_2}\in \mathcal {C}([0,1];\mR^{k+m});d(\chi_1,h)<\eta,d(\chi_2,X^h)>\alpha\}.
\de
By Lemma \ref{continuity}, we know that $h\rightarrow Z^h$ is continuous with respect to the distance $d$ when
\ce
h\in B_R=\{h\in\cH_k^2;\frac{1}{2}\|h\|_{\cH_k^2}^2\leq 2R\}.
\de
Therefore there exists $\eta>0$ such that if $\chi_1\in B_R$ and $d(\chi_1,h)\leq \eta$, then
\ce
d(X^{\chi_1},X^h)\leq\alpha.
\de
This implies that for such a value $\eta$, there is no trajectory $\chi\in \Lambda_{\eta}$ such that
\ce
\frac{1}{2}\|\chi_1\|_{\cH_k^2}^2\leq 2R\quad\text{and}\quad\chi={\chi_1\choose X^{\chi_1}}.
\de
Consequently, we have
\ce
\inf_{f\in \Lambda_{\eta}}I(f)\geq 2R,
\de
and the proof is completed by Theorem \ref{Main}.
\end{proof}

We adopt the following definition from \cite{B}.
\bd\label{contraction}
Let $U$ be an open set of $\mR^m$. For $\alpha>0$, let $\Gamma_{\alpha}:U\rightarrow U$ be a $\mathcal {C}^2$ bijective transformation having continuous derivatives up to order $2$. The family $\Gamma=\{\Gamma_{\alpha}\}_{\alpha>0}$ is said to be a system of contractions centered at $x$ if
\begin{enumerate}
\item[(a)]
$\Gamma_{\alpha}(x)=x$ for every $\alpha>0$;
\item[(b)]
if $\alpha\geq\beta$, then $|\Gamma_{\alpha}(y)-\Gamma_{\alpha}(z)|\leq|\Gamma_{\beta}(y)-\Gamma_{\beta}(z)|$ for every $y,z\in U$;
\item[(c)]
$\Gamma_1$ is the identical mapping on $U$ and $\Gamma_{\alpha^{-1}}=\Gamma_{\alpha}^{-1}$. Moreover, for every compact subset $K$ of $U$ and $\epsilon>0$, there exists $\delta>0$ such that if $|\alpha\beta-1|<\delta$ then
\ce
|\Gamma_{\alpha}\circ\Gamma_{\beta}(y)-y|<\epsilon
\de
for every $y\in K$.
\end{enumerate}
\ed

In the sequel, we require $\overline{D(A)}\subset U$.

\subsection{The large time functional iterated logarithm law}

Let now $\tilde{b}:U\rightarrow\mR^m$ and $\tilde{\sigma}:U\rightarrow\mR^{m\times k}$. Let $\widetilde{A}$ be a multivalued maximum monotone operator. We also set
\ce
&&L(u):=\log\log u,\quad\phi(u):=\sqrt{uL(u)},\quad u>e,\\
&&
\widetilde{L}=\frac{1}{2}\sum_{ij}\tilde{a}_{ij}(y)\frac{\partial^2}{\partial y_i\partial y_j}+\sum_i\tilde{b}_i(y)\frac{\partial}{\partial y_i},
\de
where $\tilde{a}=\tilde{\sigma}\tilde{\sigma}^t$.
For every $\alpha>e$, we define
\ce
&&\tilde{b}_{\alpha}(y):=\alpha(\widetilde{L}\Gamma_{\phi(\alpha)})(z)|_{z=\Gamma_{\phi(\alpha)}^{-1}(y)},\quad
\tilde{\sigma}(y):=\phi(\alpha)\nabla \Gamma_{\phi(\alpha)}(z)|_{z=\Gamma_{\phi(\alpha)}^{-1}(y)}\cdot\tilde{\sigma}(\Gamma_{\phi(\alpha)}^{-1}(y)),\\
&&\widetilde{A}_{\alpha}(y):=\alpha\nabla \Gamma_{\phi(\alpha)}(z)|_{z=\Gamma_{\phi(\alpha)}^{-1}(y)}\cdot \widetilde{A}(\Gamma_{\phi(\alpha)}^{-1}(y)).
\de

We shall make use of the following assumptions.
\begin{enumerate}[{\bf (C1)}]
\item
There exist $b:U\rightarrow\mR^m$ and $\sigma:U\rightarrow\mR^{m\times k}$  such that
\ce
\lim_{\alpha\rightarrow\infty}\tilde{b}_{\alpha}(y)=b(y),\quad
\lim_{\alpha\rightarrow\infty}\tilde{\sigma}_{\alpha}(y)=\sigma(y),
\de
uniformly on compact subsets of $U$;
\item
there exists a multivalued maximal monotone operator $A:\mR^m\rightarrow 2^{\mR^m}$ such that for every $\epsilon>0$,
\ce
\lim_{\alpha\rightarrow\infty}(1+\epsilon\widetilde{A}_{\alpha})^{-1}(y)=(1+\epsilon A)^{-1}(y)
\de
uniformly on compact subsets of $U$;
\item
$\widetilde{A}_{\alpha}$, $\widetilde{A}$, $A$ are maximal monotone operator with $\overline{D(\widetilde{A}_{\alpha})}=\overline{D(A)}$ and nonempty interior. Moreover, for $a\in \mathrm{Int}(D(\widetilde{A}_{\alpha}))$, $\widetilde{A}_{\alpha}$ is locally bounded at $a$ uniformly in $\alpha$;
\item
if $b_{\epsilon}=\tilde{b}_{1/\epsilon}$, $\sigma_{\epsilon}=\tilde{\sigma}_{1/\epsilon}$ and $A_{\epsilon}=\widetilde{A}_{1/\epsilon}$, the system of small random perturbation $(b_{\epsilon},\sigma_{\epsilon},A_{\epsilon})$ satisfies assumptions {\bf (H2)}-{\bf (H4)}.
\end{enumerate}

Consider the following MSDE:
\ce\label{origin}
\begin{cases}
\dif Y(t)\in\tilde{b}(Y(t))\dif t+\tilde{\sigma}(Y(t))\dif W(t)-\widetilde{A}(Y(t))\dif t,\\
Y(0)=x\in\overline{D(A)},
\end{cases}
\de
and let us set for $u>e$
\ce
Y^{(u)}(t)=Y(ut),\quad Z_u(t)=\Gamma_{\phi(u)}(Y^{(u)}(t)).
\de

Now we give the large time functional iterated logarithm law.
\bt\label{iterated}
Under assumptions {\bf (C1)}-{\bf (C4)}, the family $\{Z_u\}_{u>e}$ is relatively compact. Moreover, $\Theta:=\{f; I(f)\leq 1\}$ is the limit set $\{Z_u\}_{u>e}$ for $u\rightarrow\infty$ a.s., where
\ce
I(f)=\frac{1}{2}\inf_{\{h\in\cH_k^2;f=X^h\}}\|h\|^2_{\cH_k^2},
\de
and $X^h(t,x)$ solves the following MDE:
\ce
\begin{cases}
\dif X^h(t)\in b(X^h(t))\dif t+
\sigma(X^h(t))\dot h(t)\dif t-A(X^h(t))\dif t,\\
X^h(0)=x.
\end{cases}
\de
\begin{proof}
The proof of this theorem is essentially the same as the proof of Theorem 3.1 in \cite{C}. That is, one first prove that for every $c>1$ and for every $\epsilon$, there exists a.s. a positive integer $j_0=j_0(\omega)$ such that for every $j>j_0$,
\ce
d(Z_{c^j},\Theta)<\epsilon.
\de
One then applies Proposition 3.2 in \cite{C} to obtain that
\ce
\lim_{u\rightarrow\infty}d(Z_u,\Theta)=0\quad\text{a.s.}.
\de
Finally one uses Theorem \ref{exponential} to prove that for any $\epsilon>0$, there exists $c_{\epsilon}>1$ such that for every $c>c_{\epsilon}$,
\ce
\bP\{d(Z_{c^j},f)<\epsilon,\text{i.o.}\}=1,
\de
which ensures that all points in $\Theta$ are actually limit points.
Therefore we complete the proof.
\end{proof}
\et

\br
We should point out that we can not use the method developed in \cite[Proposition 2.5]{B} due to the existence of the multivalued maximal monotone operator.
\er

\subsection{The small time functional iterated logarithm law}

Without making a great effort, we can prove the functional iterated logarithm law for small value of the parameter.

Let $\tilde{b}:U\rightarrow\mR^m$ and $\tilde{\sigma}:U\rightarrow\mR^{m\times k}$. Let $\widetilde{A}$ be a multivalued maximum monotone operator. We also set
\ce
&&G(u):=\log\log u^{-1},\quad\varphi(u):=\sqrt{uG(u)},\quad 0<u<e^{-1},\\
&&\widetilde{L}:=\frac{1}{2}\sum_{ij}\tilde{a}_{ij}(y)\frac{\partial^2}{\partial y_i\partial y_j}+\sum_i\tilde{b}_i(y)\frac{\partial}{\partial y_i},
\de
where $\tilde{a}=\tilde{\sigma}\tilde{\sigma}^t$.
For every $0<\alpha<e^{-1}$, we define
\ce
\tilde{b}_{\alpha}(y)&:=&\alpha(\widetilde{L}\Gamma_{\varphi(\alpha)})(z)|_{z=\Gamma_{\varphi(\alpha)}^{-1}(y)},\\
\tilde{\sigma}_{\alpha}(y)&:=&\varphi(\alpha)\nabla \Gamma_{\varphi(\alpha)}(z)|_{z=\Gamma_{\varphi(\alpha)}^{-1}(y)}\cdot\tilde{\sigma}(\Gamma_{\varphi(\alpha)}^{-1}(y)),\\
\widetilde{A}_{\alpha}(y)&:=&\alpha\nabla \Gamma_{\varphi(\alpha)}(z)|_{z=\Gamma_{\varphi(\alpha)}^{-1}(y)}\cdot \widetilde{A}(\Gamma_{\varphi(\alpha)}^{-1}(y)).
\de

We shall make use of the following assumption.
\begin{enumerate}[{\bf (C$'$1)}]
\item
There exist $b:U\rightarrow\mR^m$ and  $\sigma:U\rightarrow\mR^{m\times k}$ such that
\ce
\lim_{\alpha\rightarrow0^+}\tilde{b}_{\alpha}(y)=b(y),\quad
\lim_{\alpha\rightarrow0^+}\tilde{\sigma}_{\alpha}(y)=\sigma(y)
\de
uniformly on compact subsets of $U$;
\item
there exists a multivalued maximal monotone operator $A:\mR^m\rightarrow 2^{\mR^m}$ such that for every $\epsilon>0$,
\ce
\lim_{\alpha\rightarrow0^+}(1+\epsilon\widetilde{A}_{\alpha})^{-1}(y)=(1+\epsilon A)^{-1}(y)
\de
uniformly on compact subsets of $U$;
\item
$\widetilde{A}_{\alpha}$, $\widetilde{A}$, $A$ are maximal monotone operator with $\overline{D(\widetilde{A}_{\alpha})}=\overline{D(A)}$ which has nonempty interior. Moreover, for $a\in \mathrm{Int}(D(\widetilde{A}_{\alpha}))$, $\widetilde{A}_{\alpha}$ is locally bounded at $a$ uniformly in $\alpha$;
\item
the system of small random perturbation $(\tilde{b}_{\alpha},\tilde{\sigma}_{\alpha},\widetilde{A}_{\alpha})$ satisfies assumptions {\bf (H2)}-{\bf (H4)}.
\end{enumerate}

Let $Y(t)$ be the solution of MSDE (\ref{origin}), and set for $0<u<e^{-1}$
\ce
Y^{(u)}(t)=Y(ut),\quad Z'_u(t)=\Gamma_{\varphi(u)}(Y^{(u)}(t)).
\de

Then we have the following theorem.
\bt\label{small time}
Under {\bf (C1)}-{\bf (C4)}, the family $\{Z_u\}_{0<u<e^{-1}}$ is relatively compact; moreover $\Theta'=\{f; I(f)\leq 1\}$ is the limit set $\{Z_u\}_{0<u<e^{-1}}$ for $u\rightarrow0^+$ a.s., where
\ce
I(f)=\frac{1}{2}\inf_{\{h\in\cH_k^2;f=X^h\}}\|h\|^2_{\cH_k^2},
\de
and $X^h(t,x)$ solves the following MDE:
\ce
\dif X^h(t)\in b(X^h(t))\dif t+
\sigma(X^h(t))\dot h(t)\dif t-A(X^h(t))\dif t,\quad
X^h(0)=x.
\de
\et

\subsection{Applications}

In what follows we will use our newly established large time functional iterated logarithm law and provide another approach to proving \cite[Theorem 3.1]{Rab}. Before stating the result, we present some related notations.

\bd
Suppose that $\cO$ is a closed convex subset of $\mR^m$, and
$I_{\cO}$ is the indicator of $\cO$, i.e.,
\ce
I_{\cO}:=
\begin{cases}
0,&\text{if}\quad x\in\cO,\\
\infty,&\text{if}\quad x\notin\cO.
\end{cases}
\de
The subdifferential of $I_\cO$ is given by
\ce
\partial I_{\cO}(x)&:=&\{y\in\mR^m;\<y,x-z\>_{\mR^m}\geq 0,\forall z\in\cO\}\\
&=&\begin{cases}
\emptyset,&\text{if}\quad x\notin\cO,\\
\{0\},&\text{if}\quad x\in \mathrm{Int}(\cO),\\
\Pi_x,&\text{if}\quad x\in\partial\cO,
\end{cases}
\de
where $\mathrm{Int}(\cO)$ is the interior of $\cO$ and $\Pi_x$ is the exterior normal cone at $x$.
\ed

It is well known  that $\partial I_\cO$ is a multivalued maximal monotone operator.

We now deal with the case of the half-space with normal reflection. Let
\ce
\mR_+^m=\{x=(x^1,x^2,\cdots,x^d);x^1\geq0\}\quad\text{and}\quad e_1=(1,0,\cdots,0).
\de
As in Anderson and Orey \cite{AO}, we define a transformation $\Gamma: \mathcal {C}([0,1];\mR^m)\rightarrow \mathcal {C}([0,1];\mR_+^m)$ as follows: for $w=(w^1,w^2,\cdots,w^m)\in \mathcal {C}([0,1];\mR^m)$,
\be\label{Gamma}
\Gamma(\omega)=(w^1+\bar{w}^1,w^2,\cdots,w^m),
\ee
where
\ce
\bar{w}^1(t)=-\inf_{s\in[0,t]}(w^1(s)\wedge0).
\de
We write $\Gamma_t(w)$ for $(\Gamma(w))_t$.

For $h\in\cH_k^2$ and $x\in \mR_+^m$, let $\widehat{X}^h(t,x)$ be the solution of the following differential equation:
\ce
\begin{cases}
\dif \widehat{X}^h(t)=b(\Gamma_t(\widehat{X}^h))\dif t+\sigma(\Gamma_t(\widehat{X}^h))\dot{h}_t\dif t,\\
\widehat{X}^h(0)=x\in \mR_+^m.
\end{cases}
\de
Define
\ce
\widehat{I}(f):=\frac{1}{2}\inf_{\{h\in\cH_k^2; f=\widehat{X}^h\}}\|h\|_{\cH_k^2}^2,
\de
with the understanding $\widehat{I}(f)=\infty$ if the above set is empty.
Now for $f\in \mathcal {C}([0,1];\mR_+^m)$, we define
\be\label{ra+}
I^+(f):=\inf\{\widehat{I}(g); g\in \mathcal {C}([0,1];\mR^m),\Gamma(g)=f\}.
\ee

Let $O=\mR_+^m$. Then $\partial O=\{x;x^1=0\}$. Denote by $\gamma$ the constant vector field defined on $\partial O$ with $\gamma(x)\equiv e_1$. Let $\tilde{b}:O\rightarrow\mR^m$ and $\tilde{\sigma}:O\rightarrow\mR^m\times \mR^k$. Consider the following Skorohod problem:
\be\label{Skorohod}
\begin{cases}
\dif Y(t)=\tilde{b}(Y(t))\dif t+\tilde{\sigma}(Y(t))\dif W(t)+1_{\partial O}(Y(t))\gamma(Y(t))\dif K(t),\\
Y(0)=x\in O,K(0)=0,
\end{cases}
\ee
where $K(t)$ is non-decreasing in $t$ and increasing only during $\Delta:=\{t;Y(t)\in\partial O\}$, $\Delta$ has Lebesgue measure zero. We denote by $\<Y(\cdot),K(\cdot)\>$ the unique solution of Eq. (\ref{Skorohod}).
As in Section 5, let us set for $u>e$
\ce
Y^{(u)}(t)=Y(ut),\quad Z_u(t)=\Gamma_{\phi(u)}(Y^{(u)}(t)).
\de

\bt\label{reflected iterated}
Under {\bf (C1)} in Section 5, the family $\{Z_u\}_{u>e}$ is relatively compact; moreover $\Theta=\{f; I^+(f)\leq 1\}$ is the limit set $\{Z_u\}_{u>e}$ for $u\rightarrow\infty$ a.s..
\et
\begin{proof}
Note that by \cite[Proposition 3.1]{C2}, $(Y(\cdot),K(\cdot))$ is the solution of the following MSDE:
\ce
\begin{cases}
\dif Y(t)\in\tilde{b}(Y(t))\dif t+\tilde{\sigma}(Y(t))\dif W(t)-\partial I_{O}(Y(t)),\\
Y(0)=x.
\end{cases}
\de

Moreover, $I^+(f)$ defined in (\ref{ra+}) is actually equivalent to $I(f)$ defined in (\ref{ra}). In fact, since
\ce
I(f)=\frac{1}{2}\inf_{\{h\in\cH_k^2;f=X^h\}}\|h\|_{\cH_k^2}^2,
\de
where $X^h(t,x)$ solves the following MDE:
\ce
\begin{cases}
\dif X^h(t)\in b(X^h(t))\dif t+
\sigma(X^h(t))\dot h(t)\dif t-\partial I_O(X^h(t))\dif t,\\
X^h(0)=x,
\end{cases}
\de
we have $
X^h(t,x)=\Gamma(\widehat{X}^h(t,x))$ (cf. \cite[Proposition 1]{AO}).
This implies $I^+(f)=I(f)$.

Now it suffices to check Assumptions {\bf (C2)}-{\bf (C3)}. Because $O=\mR_+^m$, it is obvious that the {\bf (C3)} is satisfied, i.e., $D(\widetilde{A}_{\alpha})=D(A)$, where
\ce
A(y)=\widetilde{A}(y)=I_O(y),\quad A_{\alpha}(y)=\alpha\nabla\Gamma_{\phi(\alpha)}(z)|_{z=\Gamma_{\phi(\alpha)}^{-1}}(y)\cdot I_{O}(\Gamma_{\phi(\alpha)}(y)).
\de
Next comes {\bf (C2)}. Consider the following indicator functions
\ce
I_{O}(y)\quad\text{and}\quad\alpha\nabla\Gamma_{\phi(\alpha)}(z)|_{z=\Gamma_{\phi(\alpha)}^{-1}}(y)\cdot I_{O}(\Gamma_{\phi(\alpha)}(y)).
\de
Since $O=\mR_+^m$, we have
\ce
I_{O}(y)=I_{O}(\Gamma_{\phi(\alpha)}(y))=
\begin{cases}
0,&\text{if}\quad y\in O,\\
\infty,&\text{if}\quad y\notin O,
\end{cases}
\de
and hence
\be\label{transformation1}
I_{O}(y)=\alpha\nabla\Gamma_{\phi(\alpha)}(z)|_{z=\Gamma_{\phi(\alpha)}^{-1}}(y)\cdot I_{O}(\Gamma_{\phi(\alpha)}(y))=
\begin{cases}
0,&\text{if}\quad y\in O,\\
\infty,&\text{if}\quad y\notin O.
\end{cases}
\ee
For any $\epsilon>0$, we define
\ce
I^{\epsilon}(y)=\inf_{u\in \mR^m}\{|y-u|^2/(2\epsilon)+I_O(u)\},\quad y\in O,
\de
and
\ce
I_{\alpha}^{\epsilon}(y)=\inf_{u\in \mR^m}\{|y-u|^2/(2\epsilon)+\alpha\nabla\Gamma_{\phi(\alpha)}(z)|_{z=\Gamma_{\phi(\alpha)}^{-1}}(y)\cdot I_{O}(\Gamma_{\phi(\alpha)}(u))\},\quad y\in O.
\de
Then by (\ref{transformation1}), we have
\ce
I^{\epsilon}(y)=I_{\alpha}^{\epsilon}(y)=\inf_{u\in O}\{|y-u|^2/(2\epsilon)\}.
\de
By \cite[Theorem 2.2]{Ba}, we have
\be\label{transformation2}
A^{\epsilon}=\partial I^{\epsilon}\quad\text{and}\quad A_{\alpha}^{\epsilon}=\partial I_\alpha^{\epsilon},
\ee
where $A^{\epsilon}$ and $A_{\alpha}^{\epsilon}$ are Yosida approximations of $\partial I_O(y)$ and $\partial (\alpha\nabla\Gamma_{\phi(\alpha)}(z)|_{z=\Gamma_{\phi(\alpha)}^{-1}}(y)\cdot I_{O}(\Gamma_{\phi(\alpha)}(y)))$, respectively.
Consequently, we conclude from (\ref{transformation2}) that {\bf (C2)} is satisfied.

Therefore the proof is completed by Theorem \ref{iterated}.
\end{proof}

Similarly, by applying the small time functional iterated logarithm law, we obtain the following theorem.
\bt\label{small reflected}
let us set for $0<u<e^{-1}$
\ce
Y^{(u)}(t)=Y(ut),\quad Z'_u(t)=\Gamma_{\varphi(u)}(Y^{(u)}(t)).
\de
Under {\bf (C$'$1)}, the family $\{Z_u\}_{0<u<e^{-1}}$ is relatively compact; moreover $\Theta=\{f; I^+(f)\leq 1\}$ is the limit set $\{Z_u\}_{0<u<e^{-1}}$ for $u\rightarrow0^+$ a.s..
\et

\section{Appendix: Proof of (\ref{key1})}

\bl
Assume that {\bf (H1)}-{\bf (H4)} hold. Let $X^{\epsilon,h_{\epsilon}}$ and $\widetilde{X}^{\epsilon,h_{\epsilon},\alpha}$ solve Eq. (\ref{control}) and Eq. (\ref{Yosida2}), respectively. Then for any $x\in\overline{D(A)}$, we have
\ce
\lim_{\alpha\rightarrow0}\bE\sup_{t\in[0,T]}|\widetilde{X}^{\epsilon,h_{\epsilon,\alpha}}(t)-X^{\epsilon,h_{\epsilon}}(t)|^2=0
\de
uniformly in $\epsilon>0$.
\el
\begin{proof}
We denote by $M_{\epsilon}^n(t)$ the $\mathcal{C}^{\infty}$-approximation of $M_{\epsilon}(t)$, which is defined as follows:
\ce
M_{\epsilon}^n(t):=n\int_{t-\frac{1}{n}}^{t+\frac{1}{n}}M_{\epsilon}(s-\frac{1}{n})\rho(n(t-s))\dif s,
\de
where $\rho$ is a mollifier with $supp \rho\subset(-1,1)$,
$\rho\in C^{\infty}$ and $\int_{-1}^1\rho(s)\dif s=1$. It is easy to obtain from (\ref{estimate1}) and (\ref{estimate2})
\ce
&&M_{\epsilon}^n(0)=0,\\
&&\lim_{n\rightarrow\infty}\bE\sup_{t\in[0,T]}|M_{\epsilon}^n(t)-M_{\epsilon}(t)|^2=0,\\
&&\sup_{t\in[0,T]}|M_{\epsilon}^n(t)|\leq\sup_{t\in[0,T]}|M_{\epsilon}(t)|\quad\text{a.s.},\\
&&\sup_{|t'-t|\leq\delta}|M_{\epsilon}^n(t')-M_{\epsilon}^n(t)|\leq \sup_{|t'-t|\leq\delta}|M_{\epsilon}(t')-M_{\epsilon}(t)|\quad\text{a.s.},\forall\delta>0,\\
&&\bE\sup_{t\in[0,T]}|\dot{M}_{\epsilon}^n(t)|^2+\bE\sup_{t\in[0,T]}|\ddot{M}_{\epsilon}^n(t)|^2\leq C_n,
\de
where $C_n$ is a constant which does not depend on $\epsilon$ and $\alpha$.

Let $(X^{\epsilon,h_{\epsilon},n}(\cdot,x),K^{\epsilon,h_{\epsilon},n}(\cdot,x))$ solve the following MSDE:
\ce
\begin{cases}
\dif X^{\epsilon,h_{\epsilon},n}(t)\in\dif M_{\epsilon}^n(t)-A_{\epsilon}(X^{\epsilon,h_{\epsilon},n}(t))\dif t,\\
X^{\epsilon,h_{\epsilon},n}(0)=x,
\end{cases}
\de
and $\widetilde{X}^{\epsilon,h_{\epsilon},\alpha,n}(t,x)$ solve the following differential equation:
\ce
\begin{cases}
\dif \widetilde{X}^{\epsilon,h_{\epsilon},\alpha,n}(t)=\dif M_{\epsilon}^n(t)-A_{\epsilon}^{\alpha}(\widetilde{X}^{\epsilon,h_{\epsilon},\alpha,n}(t))\dif t,\\
\widetilde{X}^{\epsilon,h_{\epsilon},\alpha,n}(0)=x.
\end{cases}
\de

Now we split our proof into three steps.

{\bf Step 1}

First we prove
\be\label{step1}
\lim_{n\rightarrow\infty}\bE\sup_{t\in[0,T]}|X^{\epsilon,h_{\epsilon},n}(t)-X^{\epsilon,h_{\epsilon}}(t)|=0
\ee
uniformly in $\epsilon$. For this, we proceed as follows. By \cite[Proposition 4.3]{C2}, we have
\ce
&&\sup_{t\in[0,T]}|X^{\epsilon,h_{\epsilon},n}(t)-X^{\epsilon,h_{\epsilon}}(t)|^2\\
&\leq&\sup_{t\in[0,T]}|M_{\epsilon}^n(t)-M_{\epsilon}(t)|^{1/2}(\sup_{t\in[0,T]}|M_{\epsilon}^n(t)-M_{\epsilon}(t)|+4|K^{\epsilon,h_{\epsilon},n}|_T^0+4|K^{\epsilon,h_{\epsilon}}|_T^0)^{1/2}.
\de
By Lemma \ref{b2}, we have
\ce
\bE|K^{\epsilon,h_{\epsilon}}|_T^0\leq C
\de
where $C$ is a constant which does not depend on $\epsilon$. We still need to verify
\ce
\bE|K^{\epsilon,h_{\epsilon},n}|_T^0\leq C,
\de
where $C$ is a constant which does not depend on $n$ and $\epsilon$. By Proposition \ref{bp} we have
\ce
&&|X^{\epsilon,h_{\epsilon},n}(t')-a|^2-|X^{\epsilon,h_{\epsilon},n}(t)-a|^2\\
&\leq&|M_{\epsilon}^n(t')+x-a|^2-|M_{\epsilon}^n(t)+x-a|^2\\
&&{}-2\gamma|K^{\epsilon,h_{\epsilon},n}|_t^{t'}+2\gamma\mu(t'-t)+2\mu\int_t^{t'}|X^{\epsilon,h_{\epsilon},n}(s)-a|\dif r\\
&&{}+2\int_t^{t'}\<M_{\epsilon}^n(s)-M_{\epsilon}(s),\dif K^{\epsilon,h_{\epsilon},n}(s)\>_{\mR^m}\\
&&{}+2\<X^{\epsilon,h_{\epsilon},n}(t')-x-M_{\epsilon}^n(t'),M_{\epsilon}^n(t')-M_{\epsilon}^n(t)\>_{\mR^m}\\
&\leq&C(1+\sup_{t\in[0,T]}|M_{\epsilon}(t)|^2+\sup_{t\in[0,T]}|X^{\epsilon,h_{\epsilon},n}(t)-a|)\\
&&{}+2(\sup_{|s'-s|\leq t'-t}|M_{\epsilon}(s')-M_{\epsilon}(s)|-\gamma)|K^{\epsilon,h_{\epsilon},n}|_{t'}^t.
\de
Hence
\ce
|X^{\epsilon,h_{\epsilon},n}(t')-a|^2-|X^{\epsilon,h_{\epsilon},n}(t)-a|^2+\frac{\gamma}{2}|K^{\epsilon,h_{\epsilon},n}|_{t'}^t\leq C(1+\sup_{t\in[0,T]}|M_{\epsilon}(t)|^2+\sup_{t\in[0,T]}|X^{\epsilon,h_{\epsilon},n}(t)-a|)
\de
on the set
\ce
\Omega_{\epsilon,n,\delta}=\{\sup_{|s'-s|\leq\delta}|M_{\epsilon}^n(s')-M_{\epsilon}^n(s)|\leq\frac{\gamma}{2}\}.
\de
Consequently, applying the same procedure in \cite[Proposition 4.9]{C2}, we obtain
\ce
|K^{\epsilon,h_{\epsilon},n}|_T^01_{\Omega_{\epsilon,n,\delta}}\leq\frac{C(1+\sup_{t\in[0,T]}|M_{\epsilon}(t)|^2)}{\delta}.
\de
For $p>6$, this yields
\ce
\bE|K^{\epsilon,h_{\epsilon},n}|_T^0&=&\sum_{k=1}^{\infty}\bE[|K^{\epsilon,h_{\epsilon},n}|_T^01_{\Omega_{\epsilon,n,\frac{1}{k+1}}\backslash\Omega_{\epsilon,n,\frac{1}{k}}}]\\
&\leq&C\sum_{k=1}^{\infty}(k+1)\bP(\Omega_{\epsilon,n,\frac{1}{k}}^c)\\
&\leq&C\sum_{k=1}^{\infty}(k+1)\bE\sup_{|t'-t|\leq\frac{1}{k}}|M_{\epsilon}(t')-M_{\epsilon}(t)|^p\\
&\leq&C\sum_{k=1}^{\infty}(k+1)(\frac{1}{k})^{p/2-1}<\infty.
\de
Therefore, we obtain the desired estimate, and hence (\ref{step1}) is satisfied.

{\bf Step 2}

Similarly to Step 1, we can prove
\be\label{step2}
\lim_{n\rightarrow\infty}\bE\sup_{t\in[0,T]}|\widetilde{X}^{\epsilon,h_{\epsilon},\alpha,n}(t)-\widetilde{X}^{\epsilon,h_{\epsilon},\alpha}(t)|=0
\ee
uniformly in $\epsilon$ and $\alpha$.

{\bf Step 3}

Let $\alpha,\beta>0$. Since  $A_{\epsilon}^{\alpha}(x)\in A_{\epsilon}(J^{\alpha}(x))$, we have
\ce
\<A_{\epsilon}^{\alpha}(x)-A_{\epsilon}^{\beta}(y),x-y\>_{\mR^m}\geq-(\alpha+\beta)\<A_{\epsilon}^{\alpha}(x),A_{\epsilon}^{\beta}(y)\>_{\mR^m},\quad\forall x,y\in\mR^m.
\de
Therefore
\ce
&&\sup_{t\in[0,T]}|\widetilde{X}^{\epsilon,h_{\epsilon},\alpha,n}(t)-\widetilde{X}^{\epsilon,h_{\epsilon},\beta,n}(t)|^2\\
&=&\sup_{t\in[0,T]}|-\int_0^t\<A_{\epsilon}^{\alpha}(\widetilde{X}^{\epsilon,h_{\epsilon},\alpha,n}(s))-A_{\epsilon}^{\beta}(\widetilde{X}^{\epsilon,h_{\epsilon},\beta,n}(s)),\widetilde{X}^{\epsilon,h_{\epsilon},\alpha,n}(s)-\widetilde{X}^{\epsilon,h_{\epsilon},\beta,n}(s)\>_{\mR^m}\dif s|\\
&\leq&(\alpha+\beta)\int_0^T|A_{\epsilon}^{\alpha}(\widetilde{X}^{\epsilon,h_{\epsilon},\alpha,n}(s)||A_{\epsilon}^{\beta}(\widetilde{X}^{\epsilon,h_{\epsilon},\beta,n}(s)|\dif s\\
&\leq&(\alpha+\beta)T\sup_{t\in[0,T]}|A_{\epsilon}^{\alpha}(\widetilde{X}^{\epsilon,h_{\epsilon},\alpha,n}(t)|\sup_{t\in[0,T]}|A_{\epsilon}^{\beta}(\widetilde{X}^{\epsilon,h_{\epsilon},\beta,n}(t)|.
\de
As in the proof of \cite[Proposition 4.7]{C2}, we have
\ce
|\frac{\dif}{\dif s}\widetilde{X}^{\epsilon,h_{\epsilon},\alpha,n}|\leq|\dot {M}_{\epsilon}^n(t)|+|A_{\epsilon}^0(x)|+\int_0^t|\ddot{M}_{\epsilon}^n(s)|\dif s,
\de
and this gives
\ce
\sup_{t\in[0,T]}|A_{\epsilon}^{\alpha}(\widetilde{X}^{\epsilon,h_{\epsilon},\alpha,n}(t))|\leq |A_{\epsilon}^0(x)|+2\sup_{t\in[0,T]}|\dot{M}_{\epsilon}^n(t)|+T\sup_{t\in[0,T]}|\ddot{M}_{\epsilon}^n(t)|.
\de
Consequently we have
\ce
\bE\sup_{t\in[0,T]}|\widetilde{X}^{\epsilon,h_{\epsilon},\alpha,n}(t)-\widetilde{X}^{\epsilon,h_{\epsilon},\beta,n}(t)|^2\leq C_n(\alpha+\beta),
\de
where $C_n$ is a constant independent of $\epsilon$, $\alpha$ and $\beta$. We know from \cite[Proposition 4.7]{C2} that
\ce
\sup_{t\in[0,T]}|\widetilde{X}^{\epsilon,h_{\epsilon},\alpha,n}(t)-{X}^{\epsilon,h_{\epsilon},n}(t)|\rightarrow 0\quad\text{a.s.}\quad\text{as}\quad \alpha\rightarrow0.
\de
Therefore we have
\be\label{step3}
\lim_{\alpha\rightarrow0}\bE\sup_{t\in[0,T]}|\widetilde{X}^{\epsilon,h_{\epsilon},\alpha,n}(t)-{X}^{\epsilon,h_{\epsilon},n}(t)|^2=0
\ee
uniformly in $\epsilon$.

Then the proof is completed by (\ref{step1})-(\ref{step3}).
\end{proof}
\vspace{1cm}

{\bf Acknowledgement}. We are very grateful to the referee for his/her careful reading of the manuscript and valuable suggestions.


\vspace{5mm}

\begin{thebibliography}{999}

\bibitem{AO}Anderson, R.F. and Orey, S.: Small random perturbation of dynamical systems with reflecting boundary. {\it Nagoya Math. J.} 60 (1976) 189-216.

\bibitem{B}Baldi, P.: Large deviations and functional iterated logarithm law for diffusion processes. {\it Probab. Theory Relat. Fields} 71 (1986) 435-453.

\bibitem{BS}Baldi, P. and Sanz, M.: Une remarque sur la th\'{e}orie des grandes d\'{e}viations. In Az\'{e}ma, J., Meyer, P.A. and Yor, M. (eds). {\it S\'{e}minaire de Probabilit\'{e}s XXV}, Lecture Notes in Math., 1485, 345-348, Springer-Verlag, Berlin, 1991.

\bibitem{BC}Baldi, P. and Chaleyat-Maurel, M.: An extension of Ventsel-Freidlin estimate, Stochastic analysis and related topic (Silivri, 1986), 305-327, Lecture Notes in Math., 1316, Springer, Berlin, 1988.

\bibitem{Ba}Barbu, V.: {\it Nonlinear Semigroups and Differential Equations in Banach Spaces}, Noord-Hoff Internet Publishing, Leyden, The Netherland, 1976.

\bibitem{BD} Barbu, V. and Da Prato, G.: The generator of the transition semigroup corresponding to a stochastic variational inequality. {\it Comm. Part. Diff. Equ.} 33 (2008) 1318-1338.

\bibitem{BD1}Budhiraja, A. and Dupuis, P.: A variational representation for positive functionals of infinite dimensional Brownian motion. {\it Probab. Math. Statist.} 20, no.1, Acta Univ. Wratislav. No. 2246 (2000), 39-61.

\bibitem{C}Caramellino, L.: Strassen's law of the iterated logarithm for diffusion processes for small time. {\it Stoch. Proc. Appl.} 74 (1998) 1-19.

\bibitem{C1}C\'epa, E.: \'Equations diff\'erentielles stochastiques multivoques. {\it S\'{e}minaire de probabilit\'{e}s XXIX}, Lecture Notes in Math., 1613, 86-107, Springer, Berlin, 1995.

\bibitem{C2}C\'epa, E.: Probleme de Skorohod multivoque. {\it Ann. of Prob.} 26 (2) (1998) 500-532.

\bibitem{CJ}C\'epa, E. and Jacquot, S.: Ergodicit\'e d'in\'egalit\'es variationnelles stochastiques. {\it Stoch. Stoch. Reports.} 63 (1997) 41-64.

\bibitem{DS}Deuschel, J.D. and Stroock, D.W.: Large Deviations. Academic Press, Boston, New York, 1988.

\bibitem{DE}Dupuis, P. and Ellis, R.S.: {\it A Weak Convergence Approach to the Theory of Large Deviations}. Wiley, New York, 1997.

\bibitem{HL}Hiriart-Urruty, J.-B. and Lemar\'echal, C.: {\it Fundamentals of Convex Analysis}, Grund. Text Ed., Springer, Berlin-Heidelberg. 2001.

\bibitem{IW}Ikeda, N. and Watanabe S.: {\it Stochastic Differential Equations and Diffusion Processes}. 2nd ed., Kodansha, Tokyo/North-Holland, Amsterdam, 1989.

\bibitem{Rab}Rabeherimanana, T.J.: Grandes d\'{e}viations et loi fonctionnelle du logarithme it\'{e}r\'{e} pour les processus de diffusions avec r\'{e}flexion. {\it Annales Math\'{e}matiques Blaise Pascal} 14 (2007) 61-76.

\bibitem{Ras}Rascanu, A.: Deterministic and stochastic differential equations in Hilbert spaces involving multivalued maximal monotone operators. {\it Panamer. Math. J.} 6 (3) (1996) 83-119.

\bibitem{RXZ}Ren, J., Xu, S. and Zhang, X.: Large deviation for multivalued stochastic differenital equations. {\it J. Theoretical Prob.} 23 (4) (2010) 1142-1156.

\bibitem{RZ1}Ren, J. and Zhang, X.: Stochastic flows for SDEs with non-Lipschitz coefficients. {\it Bull. Sci. Math.} 127 (2003) 739-754.

\bibitem{S}Skorohod, A.V.: Stochastic equations for diffusions in a bounded region. {\it Theory Prob. Appl.} 6 (1961) 264-274.

\bibitem{W}Wu, J.: Uniform large deviation for multivalued stochastic differenital equations with jumps. {\it Kyoto J. Math.} 51 (3) (2011) 535-559.

\bibitem{Z}Zalinescu, A.: Second order Hamilton-Jacobi-Bellman equations with unbounded operator. {\it Nonlinear Analysis} 75 (2012) 4784-4797.

\bibitem{Z1}Zhang, X.: Skorohod problem and multivalued stochastic evolution equations in Banach spaces. {\it Bull. Sci. Math.} 131 (2) (2007) 175-217.

\end{thebibliography}
\end{document}